\documentclass{article}
\usepackage{amsfonts,amsmath,amssymb,amscd}
\usepackage{fancybox}
\usepackage{float}
\usepackage{color}

\begin{document}

\newtheorem{definition}{Definition}[section]
\newtheorem{theorem}[definition]{Theorem}
\newtheorem{corollary}[definition]{Corollary}
\newtheorem{remark}[definition]{Remark}
\newtheorem{example}[definition]{Example}
\newtheorem{lemma}[definition]{Lemma}
\newtheorem{proposition}[definition]{Proposition}

\newcommand{\qed}{\hbox{\rule[-2pt]{3pt}{6pt}}}
\newenvironment{proof}{{\it Proof.}}{\hfill \qed}

\makeatletter
\def\tbcaption{\def\@captype{table}\caption}
\def\figcaption{\def\@captype{figure}\caption}
\makeatother

\title{Homology Groups of Symmetric Quandles and Cocycle Invariants of  
Links and Surface-links} 
\author{Seiichi Kamada\footnote{
The corresponding author.  This research is partially supported by Grant-in-Aid for Scientific Research, JSPS. } \/ and Kanako Oshiro\footnote{This research is partially supported by Grant-in-Aid for JSPS Research Fellowships for Young Scientists.}\\
Department of Mathematics, Hiroshima University, 
Hiroshima 739-8526, Japan\\
kamada@math.sci.hiroshima-u.ac.jp (S. Kamada) \\
d085317@hiroshima-u.ac.jp (K. Oshiro)}

\date{}
\maketitle

\begin{abstract}
We introduce the notion of a quandle with a good involution and its homology groups.  Carter et al. defined quandle cocycle invariants for oriented links and oriented surface-links.  By use of good involutions, 
quandle cocyle invariants can be defined for links and surface-links which are not necessarily oriented or orientable.  The invariants can be used in order to estimate the minimal triple point numbers of non-orientable surface-links. 
\end{abstract}



\section{Introduction}
The fundamental quandles are strong invariants of links, which are generalized to any codimension-two oriented manifold pairs, \cite{FR92, Joyce, Matveev}.  Carter et al. \cite{CJKLS03} defined homology groups of quandles, which are deeply related to homology groups of racks due to Fenn et al. \cite{FRS95}, cf. \cite{CJKS01b, LN03}.   Using cocycles $\theta$ of  quandles, Carter et al. \cite{CJKLS03} defined invariants of links in $3$-space or surface-links in $4$-space, called quandle cocycle invariants and denoted by $\Phi_\theta^{\rm ori}$ in this paper.   For details on quandle homologies and invariants, refer to 
\cite{AG03, CEGS05, CJKLS03, CJKS01b, FRS95, FRS07,  LN03, RS1, RS}.  
Quandle cocyle invariants can distinguish the left-handed trefoil from the right-handed trefoil (\cite{RS1, RS}), and the $2$-twist-spun trefoil from  its copy with the opposite orientation (\cite{AS05, CJKLS03, I, RS1, RS}), although the fundamental quandles themselves cannot.    In order to define quandle cocycle invariants, it is essential  that links or surface-links are oriented.   In this paper we introduce the notion of a quandle with a good involution, which is called a symmetric quandle, and its homology groups.   This notion enables us to consider quandle cocycle invariants, denoted by $\Phi_\theta$, for links or surface-links which are not necessarily oriented or orientable.  

For a given symmetric quandle $2$-cocycle $\theta$, we define a quandle cocycle invariant $\Phi_\theta$ of an unoriented link.  We show that  $\Phi_\theta$ is equal to the quandle cocycle invariant $\Phi_\theta^{\rm ori}$ of an oriented link when we give an orientation to the link arbitrarily (Theorem~\ref{thm:equality1dim}).   As a consequence, we see that $\Phi_\theta^{\rm ori}$ does not depend on the orientation of a link if $\theta$ is a quandle $2$-cocycle cohomologous to a symmetric quandle $2$-cocycle.   We can also use the fact $\Phi_\theta^{\rm ori}= \Phi_\theta$ for calculation of $\Phi_\theta^{\rm ori}$, since the calculation of $\Phi_\theta$ is sometimes simpler than that of $\Phi_\theta^{\rm ori}$ as seen in Example~\ref{example:trefoil}.  

For a given symmetric quandle $3$-cocycle $\theta$, we define a quandle cocycle invariant $\Phi_\theta$ of an unoriented or non-orientable surface-link.  We show that  when a surface-link is orientable, the invariant $\Phi_\theta$ is equal to the quandle cocycle invariant $\Phi_\theta^{\rm ori}$ of an oriented surface-link when we give any orientation to the surface-link  (Theorem~\ref{thm:equality2dim}).    The invariants $\Phi_\theta$ are essentially new and interesting for non-orientable surface-links.   As $\Phi_\theta^{\rm ori}$ can be used to estimate the minimal triple point numbers of oriented surface-links \cite{SatohShima2004, SatohShima2005}, we can use  $\Phi_\theta$ to estimate those of non-orientable surface-links 
(Theorem~\ref{theorem:triplepoint}, Proposition~\ref{proposition:cocycledihedral4}, Theorems~\ref{theorem:triple} and \ref{theorem:P2P2}).  

This paper is organized as follows:
In \S \ref{Good_involutions} the definition of a quandle with a good involution or a 
symmetric quandle is defined.
In \S \ref{Good_involutions_basic_quandles} 
we determine good involutions of trivial quandles, dihedral quandles and keis.  In \S \S \ref{The associated group of a symmetric quandle} and \ref{Homology groups of a symmetric quandle} 
the definitions of the associated group and homology groups of a symmetric quandle are given.
In \S \ref{Link_invariants_1} we investigate cocycle invariants $\Phi_\theta$ for classical links and show that 
$\Phi_\theta$ is equal to $\Phi_\theta^{\rm ori}$.  
In \S \ref{Link_invariants_2} a special case where $X$ is a kei and $\rho$ is the identity map of $X$ is considered. 
In \S \ref{Surface-link invariants} 
we discuss quandle cocycle invariants for surface-links. 
In \S \ref{Examplescocycles} examples of quandle cocycle invariants 
and applications on the minimal triple point numbers of surface-links are given. 

\section{Quandles with Good involutions}\label{Good_involutions}

A {\it quandle} \cite{FR92, Joyce,  Matveev} is a set $X$ with a binary operation $(x,y)\mapsto x^y$ satisfying the following:  (Q1) For any $x\in X $, $x^x=x$, 
(Q2) for any $x,y\in X$, there exists a unique element $z\in X$ with $z^y=x$, and 
(Q3) for any $x,y,z\in X$, $(x^y)^z=(x^z)^{(y^z)}$D
We denote the element $z$ given in (Q2) by $x^{y^{-1}}$.

A {\it rack} \cite{FR92} is a set $X$ with a binary operation satisfying (Q2) and (Q3). 

A {\it kei} \cite{Takasaki} or an {\it  involutory quandle} \cite{Joyce} is a quandle satisfying that $(x^y)^y=x$ for  any $x,y\in X$.
In other words, a kei is a quandle $X$ with $x^y=x^{y^{-1}}$ for any $x,y \in X$.

\begin{definition}{\rm 
A map $\rho :X\to X $  is {\it a good involution} \cite{ka} if it is an involution (i.e., 
$\rho \circ \rho = {\rm id}$) 
such that 
$ \rho (x^y)=\rho (x)^y $ and $ x^{\rho (y)}=x^{y^{-1}}$ for any $x, y \in X$. 
Such a pair $(X, \rho)$ is called a {\it quandle with a good involution} or a 
 {\it symmetric quandle}.  
}
\end{definition}

A ({\it quandle}) {\it homomorphism} $f:X \to Y $ is a map such that $f(x^y)=f(x)^{f(y)}$ for any $x,y\in X$.
A ({\it symmetric quandle}) {\it homomorphism} $f:(X,\rho _X)\to (Y,\rho _Y)$ means a map $f:X\to Y$ such that $f(x^y)=f(x)^{f(y)}$ and $f(\rho _X(x))=\rho _Y(f(x))$ for any $x,y\in X$. \\


\begin{example}
{\rm (cf. \cite{ka})
Let $G$ be  a group. 
The {\it conjugation quandle}, denoted by ${\rm conj}(G)$, is $G$ with the operation $x^y = y^{-1}xy$.  
The inversion, ${\rm inv}(G):G\to G; g\mapsto g^{-1}$ is a good involution of ${\rm conj}(G)$.
We call $({\rm conj}(G), {\rm inv}(G))$ the {\it conjugation symmetric quandle}.}
\end{example}

\begin{example}
\label{doublecover}
{\rm (cf. \cite{ka})
Let $X$ be a quandle and let $X_1$ and $X_2$ be two copies of $X$. For an element $x \in X$, we denote by $x_1$ and $x_2$ the corresponding elements of $X_1$ and $X_2$, respectively.  
Let $D(X)$ be the disjoint union of $X_1$ and $X_2$. For any elements $x,y \in X$, we put 
$x_i^{y_1}= (x^y)_i$ and 
$x_i^{y_2} = (x^{y^{-1}})_i$ for $i=1,2$.  
Then the set $D(X)$ is a quandle. 
The involution $\rho: D(X) \to D(X)$ interchanging  $x_1$ and  $x_2$ for every $x\in X$ is a good involution.  
}
\end{example}

\begin{definition}
{\rm (cf. \cite{ka}) 
For a codimension-two manifold pair $(W,L)$, let $\widetilde{Q}_L$ be the set of homotopy classes, $x= [(D, \alpha)]$, of all pairs $(D, \alpha)$, where $D$ is an oriented meridian disk of $L$ and $\alpha$ is a path in $W \setminus L$ starting from a point of $\partial D$ and ending at a fixed base point $\ast \in W \setminus L$.  It is a quandle with an operation defined by 
\[
[(D_1, \alpha_1)]^{[(D_2, \alpha_2)]} = [(D_1, \alpha_1 \cdot \alpha_2^{-1} \cdot 
\partial D_2 \cdot \alpha_2)].
\]
An involution $\rho:  \widetilde{Q}_L \to \widetilde{Q}_L$   defined by 
$[(D, \alpha)] \mapsto [(-D, \alpha)]$  
is a good involution of $\widetilde{Q}_L$, 
where $-D$ stands for the meridian disk $D$ with the opposite orientation. 
The {\it fundamental symmetric quandle} of $L$ is $(\widetilde{Q}_L, \rho)$. 
} 
\end{definition}

\begin{example}\label{example:unknot}{\rm 
If $L$ is an unknot in ${\mathbb R}^3$ or an unknotted $2$-sphere in ${\mathbb R}^4$,  then $\widetilde{Q}_L$ consists of two elements, say $x_1 =[(D, \alpha)]$ and $x_2=[(-D, \alpha)]$.  This is a trivial quandle (i.e., 
$x_i^{x_j}=x_i$ for $i,j \in \{1,2\}$)  and $\rho(x_1) =x_2$.  On the other hand, if $L$ is an unknotted projective plane in ${\mathbb R}^4$, then $\widetilde{Q}_L$ consists of a single element $x$ with $\rho(x)=x$, where any surface ambiently isotopic to the projective planes that project to the standard cross-cap is unknotted.
}\end{example}

When $L$ is transversely oriented (i.e., all meridian disks of $L$ are oriented coherently), we have a subquandle $Q_L$ of $\widetilde{Q}_L$ consisting of all homotopy classes of pairs $(D, \alpha)$ such that  $D$ has the given orientation; cf. p.~359 of \cite{FR92}.   Then $(\widetilde{Q}_L, \rho)$ is isomorphic to $(D(Q_L), \rho)$  in the sense of Example~\ref{doublecover}. \\


\section{Good involutions of trivial quandles, dihedral quandles and keis}\label{Good_involutions_basic_quandles}

In this section we study good involutions of trivial quandles, dihedral quandles and keis.

A {\it trivial quandle} is a quandle whose operation is trivial, i.e., $x^y=x$ for any $x,y\in X$.
A trivial quandle consisting of $n$ elements will be denoted by $T_n$ later. 

\begin{proposition}
\label{proposition:trivial}
If $X$ is a trivial quandle, then every involution of $X$ is a good involution. 
Conversely, if every involution of $X$ is a good involution, then $X$ is a trivial quandle.
\end{proposition}

\begin{proof}
Let $\rho$ be an involution of a trivial quandle $X$. For any $x, y\in X$,
$\rho (x^y)=\rho (x)=\rho (x)^y$, and 
$x^{\rho (y)}=x=x^{y^{-1}}$.
Therefore, $\rho$ is a good involution of X.

Suppose that any involution of $X$ is a good involution. 
For any $x,y\in X$, consider an involution $\rho : X\to X$ with $\rho (x)=y$, $\rho (y)=x$ and $\rho(z)=z$ for $z \not = x, y \in X$.
If $\rho $ is a good involution, then $x=\rho(y) =\rho (y^y)=\rho(y)^y=x^y$.
Thus $X$ is a trivial quandle.
\end{proof}

The {\it dihedral quandle} of order $n$ is ${\mathbb Z}/n{\mathbb Z}$ with the operation $x^y=2y-x\pmod{n}$.
We denote it by $R_n$.
For simplicity, we denote an element $i+n\mathbb{Z}$ of ${\mathbb Z}/n{\mathbb Z}$ by 
$i$ for $i \in \{0,1,\cdots ,n-1\}$. 

Operations of the dihedral quandle $R_n$ correspond to inversions along axes of symmetry of a regular $n$-gon.
Consider a regular $n$-gon whose vertices are named $0,1,\cdots ,n-1 \in \mathbb{Z}/n \mathbb{Z}=X$ in this order.
The inversion along the axis of symmetry passing through a vertex $y\in X$ induces a bijection $X \to X$ with $x\mapsto 2y-x\pmod n$.

\begin{figure}[H]
\begin{center}
\begin{minipage}{80pt}
\begin{picture}(75,75)
\put(35,35){\vector(0,1){30}}
\put(35,35){\vector(0,-1){30}}
\put(35,35){\vector(1,0){30}}
\put(35,35){\vector(-1,0){30}}
\put(35,35){\vector(1,1){21.216}}
\put(35,35){\vector(-1,1){21.216}}
\put(35,35){\vector(1,-1){21.216}}
\put(35,35){\vector(-1,-1){21.216}}

\qbezier(35,65)(35,65)(13.784,56.216)
\qbezier(5,35)(5,35)(13.784,56.216)
\qbezier(5,35)(5,35)(13.784,13.784)
\qbezier(35,5)(35,5)(13.784,13.784)
\qbezier(35,5)(35,5)(56.216,13.784)
\qbezier(65,35)(65,35)(56.216,13.784)
\qbezier(65,35)(65,35)(56.216,56.216)
\qbezier(35,65)(35,65)(56.216,56.216)

\put(33,68){$0$}
\put(5,57){$1$}
\put(-3,33){$2$}
\put(5,8){$3$}
\put(33,-5){$4$}
\put(59,8){$5$}
\put(68,33){$6$}
\put(59,57){$7$}
\end{picture}
\end{minipage}
\hspace{7mm}
\begin{minipage}{80pt}
\begin{picture}(75,75)
\put(35,35){\vector(0,1){30}}
\put(35,35){\vector(0,-1){30}}
\put(35,35){\vector(1,0){30}}
\put(35,35){\vector(-1,0){30}}

\qbezier(35,65)(35,65)(13.784,56.216)
\qbezier(5,35)(5,35)(13.784,56.216)
\qbezier(5,35)(5,35)(13.784,13.784)
\qbezier(35,5)(35,5)(13.784,13.784)
\qbezier(35,5)(35,5)(56.216,13.784)
\qbezier(65,35)(65,35)(56.216,13.784)
\qbezier(65,35)(65,35)(56.216,56.216)
\qbezier(35,65)(35,65)(56.216,56.216)

\put(33,68){$0$}
\put(5,57){$1$}
\put(-3,33){$2$}
\put(5,8){$3$}
\put(33,-5){$4$}
\put(59,8){$5$}
\put(68,33){$6$}
\put(59,57){$7$}
\end{picture}
\end{minipage}
\hspace{7mm}
\begin{minipage}{80pt}
\begin{picture}(75,75)
\put(35,35){\vector(1,1){21.216}}
\put(35,35){\vector(-1,1){21.216}}
\put(35,35){\vector(1,-1){21.216}}
\put(35,35){\vector(-1,-1){21.216}}

\qbezier(35,65)(35,65)(13.784,56.216)
\qbezier(5,35)(5,35)(13.784,56.216)
\qbezier(5,35)(5,35)(13.784,13.784)
\qbezier(35,5)(35,5)(13.784,13.784)
\qbezier(35,5)(35,5)(56.216,13.784)
\qbezier(65,35)(65,35)(56.216,13.784)
\qbezier(65,35)(65,35)(56.216,56.216)
\qbezier(35,65)(35,65)(56.216,56.216)

\put(33,68){$0$}
\put(5,57){$1$}
\put(-3,33){$2$}
\put(5,8){$3$}
\put(33,-5){$4$}
\put(59,8){$5$}
\put(68,33){$6$}
\put(59,57){$7$}
\end{picture}
\end{minipage}
\end{center}

\figcaption{The antipodal map and the two half-antipodal maps}
\label{g.i.d}
\end{figure}
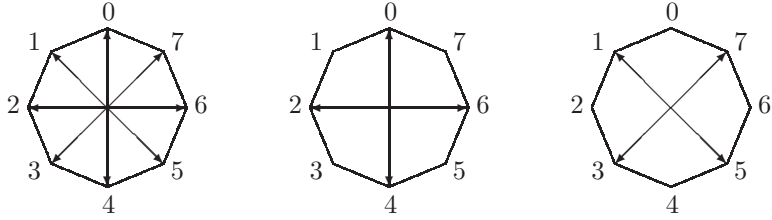

Let $X$ be the dihedral quandle $R_n=\mathbb{Z}/n\mathbb{Z}$.
The identity map is a good involution.
If $n$ is an even number, say $2m$, then a map $\rho :X\to X$ with $\rho(i)=i+m$ for all $i$ is a good involution.
We call it the {\it antipodal map}. 
Moreover if $n=2m$ for some even number $m$, then $X$ has two more good involutions,
say $\rho_1$ and $\rho_2$, with $\rho_1(i)=i$ for odd $i$'s and $\rho_1(i)=i+m$ for even $i$'s;
and $\rho_2(i)=i+m$ for odd $i$'s and $\rho_2(i)=i$ for even $i$'s. 
We call them {\it half-antipodal maps}. See Figure~\ref{g.i.d}.  

\begin{theorem}
\label{theorem:dihedral}
Let $X$ be the dihedral quandle $R_n$.

{\rm (1)} When $n$ is an odd number, the identity map is the only good involution.

{\rm (2)} When $n=2m$ and $m$ is an odd number, a good involution is either the identity map or the antipodal map.

{\rm (3)} When $n=2m$ and $m$ is an even number, there are four good involutions: the identity map, the antipodal map and two half-antipodal maps.
\end{theorem}

We prepare a lemma for the proof of Theorem \ref{theorem:dihedral}.
\begin{lemma}
\label{sec:lemma1}
Let $\rho $ be a good involution of $R_n$. 

{\rm (1)} For each $i\in R_n$, $\rho (i)=i$ or $\rho (i)=i+n/2$. 
The latter occurs only if $n$ is even. 

{\rm (2)} 
Let $n$ be an even number.
If  $\rho (i)=i+n/2$ for some $i$, then  
 $\rho (j)=j+n/2$ for every $j$ such that $ i \equiv j \pmod 2$. 
\end{lemma}

\begin{proof} (1)  
For any $i$ and $j\in R_n$, $j^{\rho(i)}=j^{i^{-1}}=j^i$.
Thus $2 \rho(i)=2i$ in ${\mathbb Z}/n{\mathbb Z}$. If $n$ is odd, then $\rho(i)=i$.
If $n$ is even, then $\rho(i)=i$ or $i+n/2$.

(2) Suppose  that $\rho (i)=i+n/2$ for some $i$.  Then 
$\rho(i+2)=\rho(i^{i+1})=\rho(i)^{i+1}=(i+n/2)^{i+1}=i+2-n/2=i+2+n/2$.  By induction, we have the result. 
\end{proof}

\vspace{10pt} 

{\it Proof of Theorem~{\rm \ref{theorem:dihedral}}.}
Let $\rho $ be a good involution of $R_n$.

(1) When $n$ is an odd number, 
by Lemma~\ref{sec:lemma1} (1), $\rho $ is the identity map. 

(2) When $n=2m$ for an odd number $m$, suppose that $\rho$ is not the identity map.
By Lemma~\ref{sec:lemma1} (1), there is an element $i\in R_n$ such that $\rho(i)=i+m$.   
Then, $\rho (i+m)=i$. 
One of $i$ and $i+m$ is an odd number, and the other is an even number.
By Lemma~\ref{sec:lemma1} (2), $\rho (j)=j+m$ for every $j\in R_n$.

(3) When $n=2m$ for an even number $m$, suppose that $\rho$ is not the identity map.
By Lemma~\ref{sec:lemma1} (1), there is an element $i \in R_n$ with $\rho(i)=i+m$. 
By  Lemma~\ref{sec:lemma1} (2), we see that 
$\rho$ is the antipodal or a half-antipodal map.
\hfill\qed

A good involution of a quandle is not necessarily a quandle homomorphism.  It is a quandle homomorphism if and only if the quandle is a kei.  A kei is characterized in terms of a good involution as follows.
\begin{proposition}
\label{kei-goodinvolution}
Let $X$ be a quandle.  The following are equivalent.  

{\rm (1)} $X$ is a kei. 

{\rm (2)} The identity map of $X$ is a good involution.

{\rm (3)} $X$ has a good involution which is a quandle homomorphism. 

{\rm (4)} Any good involution of $X$ is a quandle homomorphism, and $X$ has at least one good involution.
\end{proposition}

\begin{proof}
(1) $\Rightarrow$ (2):  For the identity map $\rho$, $\rho(x^y) = \rho(x)^y$.  Since $X$ is a kei, 
$x^{\rho(y)}=x^y=x^{y^{-1}}$.

(2) $\Rightarrow$ (3): It is obvious. 

(3) $\Rightarrow$ (1): Let $\rho$ be a good involution which is a quandle homomorphism.
For any $x,y \in X$, $\rho (x^y)=\rho(x)^{\rho (y)}$.
Since $\rho(x^y)=\rho (x)^y$ and $\rho (x)^{\rho (y)}=\rho (x)^{y^{-1}}$,
we have $\rho (x)^y=\rho(x)^{y^{-1}}$.

(1) $\Rightarrow$ (4):  Let $\rho$ be a good involution. 
For any $x,y\in X$,
$\rho (x^y)=\rho(x)^y=\rho(x)^{y^{-1}}=\rho(x)^{\rho(y)}.$
Therefore, $\rho$ is a quandle homomorphism.  As seen above, the identity map is a good involution.

(4) $\Rightarrow$ (3): It is obvious. \end{proof}


\section{The associated group of a symmetric quandle}\label{The associated group of a symmetric quandle}

The {\em associated group} of a quandle $X$  is  
\[
G_X = \langle x \in X \/ ; \/ 
x^y = y^{-1} x y  \quad (x, y \in X) \/ \rangle.  
\]
The associated group with the natural map $\eta: X \to G_X$ has a universal mapping property; see \cite{FR92, Joyce}.  

The {\it associated group}, $G_{(X, \rho)}$, of a symmetric quandle $(X,\rho)$ is defined by 
\[
G_{(X, \rho)} = \langle x \in X \/ ; \/ 
x^y = y^{-1} x y  \quad (x, y \in X),   \quad 
\rho(x) = x^{-1}  \quad (x \in X)   \/ 
 \rangle;  
\]
The {\it natural map} $\mu: X \to G_{(X, \rho)}$ is the composition of the inclusion map 
$X \to F(X)$ and the projection map $F(X) \to G_{(X, \rho)}$, where 
$F(X)$ is the free group on $X$.  
The associated group $G_{(X, \rho)}$ with
the natural map  $\mu: X \to G_{(X, \rho)}$ has a universal mapping property as follows: 

\begin{proposition}
\label{sec:conjugation}
Let $(X,\rho)$ be a symmetric quandle.

{\rm (1)} The natural map $\mu : X\to G_{(X, \rho)}$ is a symmetric quandle homomorphism $$\mu: (X,\rho ) \to \big({\rm conj}(G_{(X, \rho)}), {\rm inv}(G_{(X, \rho)})\big).$$

{\rm (2)} Let $G$ be a group.
Any symmetric quandle homomorphism $f: (X,\rho)\to ({\rm conj}(G),{\rm inv}(G))$ ``factors uniquely through $\mu$''; that is, 
there exists a unique group homomorphism $f_\sharp: {G_{(X, \rho)}} \to G$ with $f=f_\sharp \circ \mu$, i.e., it makes the following diagram commutative.
\begin{center}
$$
\begin{CD}
(X,\rho) @>\mu >> G_{(X, \rho)}\\
@VfVV @VVf_\sharp V \\
({\rm conj}(G), {\rm inv}(G)) @>{\rm id} >>G
\end{CD}
$$
\end{center}
\end{proposition}

\begin{proof}
(1) For any $x, y \in X$, $x^y = y^{-1}xy$ and $\rho(x) = x^{-1}$ in $G_{(X, \rho)}$.  This implies that $\mu(x^y) = \mu(x)^{\mu(y)}$ and $\mu \circ \rho(x) = {\rm inv}(G_{(X, \rho)}) \circ \mu(x)$.  

(2) Let $f_\ast :  F(X) \to G$ be the group homomorphism determined by 
the map $f:X\to G$.  
Since $f$ is a symmetric quandle homomorphism, 
$$
f(x^y)=f(y)^{-1}f(x)f(y)\quad \mbox{and}\quad f(\rho(x))=f(x)^{-1}
$$
for all $x,y\in X$.  Thus $f_\ast :F(X) \to G$ induces a desired homomorphism 
$f_\sharp :G_{(X, \rho)}\to G$.
\end{proof}

Proposition~\ref{sec:conjugation} is a symmetric quandle version of Proposition 2.1 of \cite{FR92}. 

\begin{example}{\rm 
(1) Let $(X, \rho)$ be the fundamental symmetric quandle of 
an unknot in ${\mathbb R}^3$ or an unknotted $2$-sphere in ${\mathbb R}^4$, i.e., $X = T_2=\{ x_1, x_2\}$ and $\rho(x_1) =x_2$ (Example~\ref{example:unknot}).  The associated group $G_X$ of $X$ is ${\mathbb Z}\times {\mathbb Z}$ and the associated group 
$G_{(X, \rho)}$ is ${\mathbb Z}$.  

(2) Let $(X, \rho)$ be the fundamental symmetric quandle of an unknotted projective plane in ${\mathbb R}^4$, i.e., $X= \{ x\}$ and $\rho(x)=x$.  The associated group  $G_X$ is ${\mathbb Z}$ and the associated group 
$G_{(X, \rho)}$ is ${\mathbb Z}/2{\mathbb Z}$.  
}\end{example}

For a quandle $X$, an {\em $X$-set} is a set $Y$ equipped with a right action of the associated group 
$G_X$.   For a symmetric quandle $(X, \rho)$, 
an {\em $(X, \rho)$-set} is a set $Y$ equipped with a right action of the associated group 
$G_{(X, \rho)}$.  
We denote by $y^g$ or by $y\cdot g$ the image of an element $y\in Y$ by the action $g \in G_{(X, \rho)}$.
Then $y\cdot (x_1 x_2)=(y\cdot x_1)\cdot x_2$, $y\cdot (x_1^{x_2})=y\cdot (x_2^{-1} x_1 x_2)$ and $y\cdot (\rho (x_1))=y\cdot (x_1^{-1})$ for $x_1,x_2 \in X$ and $y\in Y$.

\section{Homology groups of a symmetric quandle}\label{Homology groups of a symmetric quandle}

Let $(X, \rho)$ be a symmetric quandle.  
Let $C_n= C_n(X)$ be the free abelian group generated by $n$-tuples $(x_1,\cdots ,x_n)$ of elements of $X$ when $n$ is a positive integer, and let $C_n$ be $\{ 0\}$ otherwise. 
Define the boundary homomorphism $\partial _n:C_n\to C_{n-1}$ by
$$
\partial _n (x_1,\cdots ,x_n) = \sum_{i=1}^n \ (-1)^i 
\{ (x_1,\cdots ,\widehat{x_i},\cdots ,x_n)
-(x_1^{x_i},\cdots ,x_{i-1}^{x_i},\widehat{x_i},x_{i+1},\cdots ,x_n) \}
$$
for $n> 1$ and $\partial _n = 0$ for $n\leq 1$. 
Then $C_\ast = \{ C_n,\partial _n \}$ is a chain complex (cf. \cite{CJKLS03, CJKS01b, FRS95, FRS07}).

Let $D_n^{\rm Q}$ be the subgroup of $C_n$ generated by the elements of
\[
\bigcup_{i=1}^{n-1} \{ (x_1,\cdots , x_n)\in X^n  \  |\    x_i=x_{i+1}  \}, 
\]
and let  $D_n^{\rho }$ be the subgroup of $C_n$ generated by the elements of 
\[
\bigcup_{i=1}^{n} \left\{
(x_1,\cdots ,x_n)+ 
(x_1^{x_i},\cdots , x_{i-1}^{x_i} ,\rho (x_i) ,x_{i+1},\cdots ,x_n) \ | \ 
 x_1,\cdots ,x_n \in X 
\right\}. 
\]
Then $D_\ast ^{\rm Q} = \{ D_n^{\rm Q},\partial _n \}$ and $D_\ast ^{\rho} = \{ D_n^{\rho },\partial _n  \}$ are subcomplexes of $C_\ast$.  The former of this fact is proved in \cite{CJKLS03} and the latter is seen as a special case of Lemma~\ref{lemma:subcomplex}.

Define abelian groups $C_n^{\rm R}(X)$, $C_n^{\rm Q}(X)$, $C_n^{{\rm R},\rho }(X)$ and $C_n^{{\rm Q},\rho }(X)$ by 
\[
\begin{array}{llllll}
C_n^{\rm R}(X)\ &=\ &C_n,  &C_n^{\rm Q}(X)\ &=\ &C_n/D_n^{\rm Q},\\
C_n^{{\rm R}, \rho }(X)\ &=\ &C_n/D_n^\rho, ~ {\rm and} & 
C_n^{{\rm Q}, \rho }(X)\ &=\ &C_n/( D_n^{\rm Q}+D_n^\rho  ), 
\end{array}
\]
and we have chain complexes $C_\ast ^{\rm R}(X)$, $C_\ast ^{\rm Q}(X)$, $C_\ast ^{{\rm R},\rho }(X)$, and $C_\ast ^{{\rm Q},\rho }(X)$.
Their homology groups are denoted by 
$H_*^{\rm R}(X)$, $H_*^{\rm Q}(X)$, $H_*^{{\rm R},\rho }(X)$, and $H_*^{{\rm Q},\rho }(X)$, and called the {\it rack homology groups}, {\it quandle homology groups}, {\it symmetric rack homology groups}, and {\it symmetric quandle homology groups}, respectively.

Let $Y$ be an $(X, \rho)$-set.   Let $C_n(X)_Y$ be the free abelian group generated by the elements
$(y, x_1, \dots, x_n)$ where $y \in Y$ and $x_1, \dots, x_n \in X$ when $n$ is a positive integer, $C_0(X)_Y={\mathbb Z}(Y)$, the free abelian group on $Y$, and let $C_n(X)_Y$ be $\{ 0\}$ for $n<0$. 

Define the boundary homomorphism $\partial_n : C_n(X)_Y \to C_{n-1}(X)_Y$ by 
$$
\partial _n (y, x_1, \cdots ,x_n) = \sum_{i=1}^n (-1)^i 
\{ (y, x_1,\cdots ,\widehat{x_i},\cdots ,x_n) 
-(y^{x_i}, x_1^{x_i},\cdots ,x_{i-1}^{x_i},\widehat{x_i},x_{i+1},\cdots ,x_n)\}
$$
for $n \geq  1$, and $\partial_n=0$ otherwise.  
Then $C_*(X)_Y= \{ C_n(X)_Y, \partial_n \}$ is a chain complex (cf. \cite{FRS95,FRS07}).  

Let $D_n^{\rm Q}(X)_Y$ be the subgroup of $C_n(X)_Y$ 
generated by the elements of
\[ 
\bigcup_{i=1}^{n-1} \{ (y, x_1, \dots, x_n)\ | \ x_i= x_{i+1}  \} 
\]
and let $D_n^{\rho}(X)_Y$ be the subgroup of $C_n(X)_Y$ 
generated by the elements of
\[
\bigcup_{i=1}^{n} 
\big\{ (y, x_1, \dots, x_n) + (y^{x_i}, x_1^{x_i}, \dots, x_{i-1}^{x_i}, \rho(x_i), x_{i+1}, \dots, x_n) \ |\ y\in Y,x_1,\cdots ,x_n \in X \big\}.
\]

\begin{lemma}
\label{lemma:subcomplex}
For each $n$, $\partial_n( D_n^{\rm Q}(X)_Y) \subset D_{n-1}^{\rm Q}(X)_Y$ and 
$\partial_n( D_n^{\rho}(X)_Y) \subset D_{n-1}^{\rho}(X)_Y$.  
Hence $D_*^{\rm Q}(X)_Y=\{ D_n^{\rm Q}(X)_Y,\partial _n \}$ and 
$D_*^{\rho}(X)_Y=\{ D_n^{\rho}(X)_Y,\partial _n \}$ are subcomplexes of $C_*(X)_Y$. 
\end{lemma}

\begin{proof}
 It is seen by a direct calculation that $\partial_n( D_n^{\rm Q}(X)_Y) \subset D_{n-1}^{\rm Q}(X)_Y$ (cf. \cite{CJKLS03}).   We show that  
$\partial_n( D_n^{\rho}(X)_Y) \subset D_{n-1}^{\rho}(X)_Y$.  
Let $i\in \{ 1,\cdots ,n \}$ be fixed. Recall that $(x^z)^{(y^z)}={(x^y)}^z$ for $x\in X\cup Y $ and $y,z\in X$, 
and we denote $(x^y)^z$ by $x^{y z}$ for simplicity.
In the following calculation, $\equiv $ stands for the congruence modulo $D_{n-1}^{\rho}(X)_Y$.

\begin{align*}
\partial _n(y^{x_i}&,x_1^{x_i},\cdots , x_{i-1}^{x_i} , \rho (x_i) ,x_{i+1},\cdots ,x_n)\vspace{2mm} \\
&\begin{array}{rl}
=\displaystyle\sum_{j=1}^{i-1} (-1)^j 
&\{ (y^{x_i},x_1^{x_i},\cdots ,\widehat{x_j^{x_i}} ,\cdots ,x_{i-1}^{x_i},\rho (x_i) ,x_{i+1},\cdots ,x_n)\vspace{2mm} \\
& -(y^{x_j x_i},x_1^{x_j x_i},\cdots ,x_{j-1}^{x_j x_i},\widehat{x_j^{x_i}},x_{j+1}^{x_i},\cdots ,x_{i-1}^{x_i},\rho (x_i),x_{i+1},\cdots ,x_n)\}\vspace{2mm} \\
+(-1)^i 
&\{ ( y^{x_i},x_1^{x_i},\cdots ,x_{i-1}^{x_i},x_{i+1}, \cdots ,x_n)
-(y,x_1,\cdots ,x_{i-1},x_{i+1}, \cdots ,x_n)\} \vspace{2mm}\\
+\displaystyle\sum_{j=i+1}^{n} (-1)^j 
&\{ (y^{x_i},x_1^{x_i},\cdots ,x_{i-1}^{x_i},\rho (x_i),x_{i+1}, \cdots ,\widehat{x_j},\cdots,x_n)\vspace{2mm} \\
&-(y^{x_i {x_j}},x_1^{x_i {x_j}},\cdots ,x_{i-1}^{x_i {x_j}},\rho (x_i^{x_j}),x_{i+1}^{x_j}, \cdots ,x_{j-1}^{x_j},\widehat{x_j},x_{j+1}\cdots,x_n)\}
\end{array}\\
&\begin{array}{rl}
\equiv \displaystyle\sum_{j=1}^{i-1} (-1)^j
&\{-(y,x_1,\cdots ,\widehat{x_j},\cdots ,x_n)+(y^{x_j},x_1^{x_j},\cdots ,x_{j-1}^{x_j},\widehat{x_j},x_{j+1},\cdots ,x_n)\}\vspace{2mm} \\
+(-1)^i&\{ ( y^{x_i},x_1^{x_i},\cdots ,x_{i-1}^{x_i},x_{i+1}, \cdots ,x_n)
-(y,x_1,\cdots ,x_{i-1},x_{i+1}, \cdots ,x_n)\} \vspace{2mm}\\
+\displaystyle\sum_{j=i+1}^{n} (-1)^j 
&\{ -(y,x_1,\cdots ,\widehat{x_j},\cdots ,x_n)
+(y^{x_j},x_1^{x_j},\cdots ,x_{j-1}^{x_j},\widehat{x_j},x_{j+1},\cdots ,x_n)\}
\end{array}\\
&\equiv -\partial _n(y,x_1,\cdots ,x_n).
\end{align*}\\
Therefore, $\partial_n( D_n^{\rho}(X)_Y) \subset D_{n-1}^{\rho}(X)_Y$. 
\end{proof}

Define $C_n^{\rm R}(X)_Y$, $C_n^{\rm Q}(X)_Y$, 
$C_n^{{\rm R}, \rho}(X)_Y$,  and $C_n^{{\rm Q}, \rho}(X)_Y$ by 
$$ \begin{array}{llllll}
C_n^{\rm R}(X)_Y &=&  C_n(X)_Y, & 
C_n^{\rm Q}(X)_Y &=&  C_n(X)_Y / D_n^{\rm Q}(X)_Y,\\
C_n^{{\rm R}, \rho}(X)_Y &=&  C_n(X)_Y / D_n^{\rho}(X)_Y, &  
C_n^{{\rm Q}, \rho}(X)_Y &=&  C_n(X)_Y / (D_n^{\rm Q}(X)_Y + D_n^{\rho}(X)_Y), 
\end{array} $$
and we have chain complexes $C_\ast^{\rm R}(X)_Y$, 
$C_\ast^{\rm Q}(X)_Y$,
$C_\ast^{{\rm R}, \rho}(X)_Y$ and 
$C_\ast^{{\rm Q}, \rho}(X)_Y$. 
The homology groups are denoted by 
$H_*^{\rm R}(X)_Y$, $H_*^{\rm Q}(X)_Y$, $H_*^{{\rm R}, \rho}(X)_Y$, and $H_*^{{\rm Q}, \rho}(X)_Y$, respectively.  

For an abelian group $A$, we define the chain and cochain complexes 
$C_*^W (X, A)_Y = C_*^W(X)_Y \otimes A$ and 
$C^*_W(X,A)_Y = {\rm Hom} (C_*^W(X)_Y, A)$, where 
$W={\rm R}, {\rm Q}, \{ {\rm R}, \rho\}$ and $\{ {\rm Q}, \rho\}$. 
  The homology and cohomology groups are denoted by 
 $H_*^W(X, A)_Y$ and $H^*_W(X, A)_Y$, respectively. 


\section{Quandle cocycle invariants $\Phi_\theta$ of classical links}\label{Link_invariants_1}

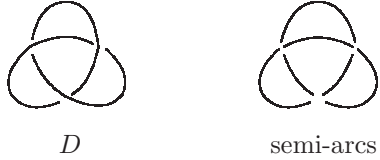
\begin{figure}[H]
\begin{center}
\begin{minipage}{60pt}
\begin{picture}(60,70)
 \qbezier(42,48)(40,65)(30,65)
 \qbezier(17,51)(20,65)(30,65)
 \qbezier(16,48)(5,38)(9,30)
 \qbezier(9,30)(15,22)(27,27)
 \qbezier(45,47)(55,38)(51,30)
 \qbezier(51,30)(45,22)(30,29)

 \qbezier(40,49)(30,55)(16,48)
 \qbezier(17,45)(21,34)(30,29)
 \qbezier(32,30)(39,34)(42,48)

\put(27,8){$D$}
\end{picture}
\end{minipage}
\hspace{10mm}
\begin{minipage}{60pt}
\begin{picture}(60,70)
 \qbezier(44,50)(40,65)(30,65)
 \qbezier(16,50)(20,65)(30,65)
 \qbezier(15,47)(5,38)(9,30)
 \qbezier(9,30)(15,22)(27,27)
 \qbezier(45,47)(55,38)(51,30)
 \qbezier(51,30)(45,22)(33,27)

 \qbezier(43,48)(30,55)(17,48)
 \qbezier(17,45)(21,34)(28,30)
 \qbezier(32,30)(39,34)(43,45)

\put(12,8){semi-arcs}
\end{picture}
\end{minipage}
\caption{semi-arcs}
\label{fgsemiarcs}
\end{center}
\end{figure}

Let $D$ be a diagram in ${\mathbb R}^2$ of an unoriented link in ${\mathbb R}^3$. 
Divide over-arcs at the crossings and we call the arcs of the result {\em semi-arcs} of $D$. 
For example, the diagram $D$ in Figure~\ref{fgsemiarcs} consists of $3$ over-arcs and $6$ semi-arcs.

\begin{figure}[H]
\begin{center}
\begin{minipage}{110pt}
\begin{picture}(110,60)
\put(0,30){\line(1,0){25}}
\put(35,30){\line(1,0){25}}
\put(30,0){\line(0,1){60}} 
\put(20,15){\vector(1,0){20}} 
\put(20,45){\vector(1,0){20}} 
\put(45,13){$x_1$} 
\put(45,43){$x_2$} 
\put(65,20){\framebox(40,14){$x_1 = x_2$} }
\end{picture}
\end{minipage}
\hspace{10mm}
\begin{minipage}{110pt}
\begin{picture}(110,60)
\put(0,30){\line(1,0){25}}
\put(35,30){\line(1,0){25}}
\put(30,0){\line(0,1){60}} 
\put(20,15){\vector(1,0){20}} 
\put(40,45){\vector(-1,0){20}} 
\put(45,13){$x_1$} 
\put(45,43){$x_2$} 
\put(65,20){\framebox(50,14){$x_1 = \rho(x_2)$} }
\end{picture}
\end{minipage}

\vspace{3mm}
\begin{minipage}{110pt}
\begin{picture}(110,60)
\put(0,30){\line(1,0){25}}
\put(35,30){\line(1,0){25}}
\put(30,0){\line(0,1){60}} 
\put(20,15){\vector(1,0){20}} 
\put(10,20){\vector(0,1){20}}
\put(50,20){\vector(0,1){20}}
\put(35,5){$x_3$} 
\put(45,43){$x_2$} 
\put(15,43){$x_1$}
\put(65,20){\framebox(40,14){$x_1^{x_3} = x_2$}}
\end{picture}
\end{minipage}
\hspace{10mm}
\begin{minipage}{110pt}
\begin{picture}(110,60)
\put(0,30){\line(1,0){25}}
\put(35,30){\line(1,0){25}}
\put(30,0){\line(0,1){60}} 
\put(20,15){\vector(1,0){20}} 
\put(10,20){\vector(0,1){20}}
\put(50,40){\vector(0,-1){20}}
\put(35,5){$x_3$} 
\put(45,43){$x_2$} 
\put(15,43){$x_1$}
\put(65,20){\framebox(55,14){$x_1^{x_3} = \rho(x_2)$}}
\end{picture}
\end{minipage}
\end{center}
\caption{Coloring conditions}
\label{fgadmcolora}
\end{figure}
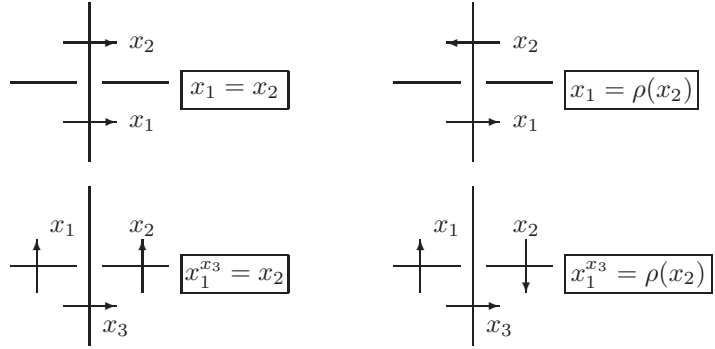

We say that an assignment of a normal orientation and an element of $X$ to each semi-arc of $D$ 
satisfies the {\em coloring conditions} if it satisfies the following two conditions  (Figure~\ref{fgadmcolora}): 
\begin{itemize}
\item 
Suppose that two diagonal semi-arcs coming from an over-arc of $D$ at a crossing $v$ are labeled by $x_1$ and $x_2$. If the normal orientations are coherent, then $x_1=x_2$, otherwise $x_1= \rho(x_2)$. 
\item 
Suppose that two diagonal semi-arcs $e_1$ and $e_2$ which are under-arcs at a crossing $v$ are labeled by $x_1$ and $x_2$, and suppose that one of the semi-arcs coming from an over-arc of $D$ at $v$, say $e_3$, is labeled by $x_3$.  We assume that the normal orientation of the over semi-arc $e_3$ points from $e_1$ to $e_2$.  If the normal orientations of $e_1$ and $e_2$ are coherent, then $x_1^{x_3}= x_2$, otherwise $x_1^{x_3}= \rho(x_2)$. 
\end{itemize}

A {\em  basic inversion} is an operation which reverses the normal orientation of a semi-arc and changes the element $x$ assigned the arc by $\rho(x)$.  Note that the coloring conditions are preserved under basic inversions.

An {\em $(X, \rho)$-coloring} of $D$ is the equivalence class of 
an assignment of a normal orientation and an element of $X$ to each semi-arc of $D$ satisfying the coloring conditions. 
Here the equivalence relation is generated by basic inversions.  

\begin{figure}[H]
\begin{center}
\begin{picture}(160,30)(0,0)
\put(0,0){\line(1,0){60}} 
\put(100,0){\line(1,0){60}} 
\put(30,10){\vector(0,-1){20}} 
\put(130,-10){\vector(0,1){20}} 
\put(77,-3){=} 
\put(35,-13){$x$} 
\put(135,-13){$\rho (x)$}
\end{picture}
\end{center}
\caption{A basic inversion}
\label{fgbasicinversion}
\end{figure}
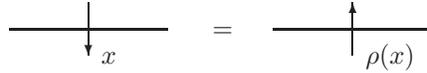

\begin{figure}[H]
\begin{center}
\begin{picture}(110,30)(0,0)
\put(0,10){\line(1,0){60}} 
\put(30,20){\vector(0,-1){20}} 
\put(22,14){$x$} 
\put(40,20){$[y_1]$}
\put(40,-4){$[y_2]$}
\put(75,3){\framebox(40,14){$y_1^{x} = y_2$} }
\end{picture}
\end{center}
\caption{Coloring condition for regions}
\label{fgadmcolorb}
\end{figure}
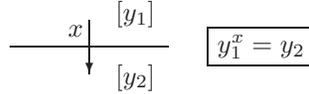

Let $Y$ be an $(X, \rho)$-set. An {\em $(X, \rho)_Y$-coloring} of $D$ is an 
$(X, \rho)$-coloring with an assignment of an element of $Y$ to each complementary region of $D$ 
satisfying the following condition  (Figure~\ref{fgadmcolorb}): 
\begin{itemize}
\item Suppose that two adjacent regions $f_1$ and $f_2$ which are separated by a semi-arc, say $e$, are labeled by $y_1$ and $y_2$.  Suppose that the semi-arc $e$ is labeled by $x$.  If the normal orientation of $e$ points from $f_1$ to $f_2$, then $y_1^x = y_2$. 
\end{itemize}

In figures, the label $y \in Y$ assigned a region is indicated by $[y]$. See Figure~\ref{fgadmcolorb}. 

Considering an $(X, \rho)$-set $Y$ consisting of a single element, we may regard an $(X, \rho)$-coloring as an $(X, \rho)_Y$-coloring. 

\begin{proposition} 
\label{proposition:coloring}
Let $(X,\rho)$ be a symmetric quandle, and $Y$ an $(X,\rho)$-set. 
If two link diagrams represent the same unoriented link type, then there is a bijection between the sets of  $(X,\rho)$-colorings of the diagrams, and there is a bijection between the sets of $(X, \rho)_Y$-colorings.
\end{proposition}

\begin{proof}
Suppose that $D$ and $D'$ are diagrams related by a single Reidemeister move.
Let $E$ be a $2$-disk in $\mathbb{R}^2$ in which the Reidemeister move is applied.
For each $(X,\rho)_Y$-coloring of $D$, its restriction to $D\setminus E (=D'\setminus E)$ 
can be uniquely extended to an $(X,\rho)_Y$-coloring of $D'$.
Thus there is a bijection between the set of $(X,\rho)_Y$-colorings of $D$ and that of $D'$. 
\end{proof}

Let $D$ be an unoriented link diagram.  
Fix an $(X, \rho)_Y$-coloring of $D$, say $C$. 
For a crossing $v$ of $D$, there are four complementary regions of $D$ around $v$. 
(Some of them may be the same.)
Choose one of them, say $f$, which we call {\it a specified region} for $v$, and let $y$ be the label of $f$.  
Let $e_1$ and $e_2$ be the under semi-arc and the over semi-arc at $v$, respectively,  which face the region $f$.  
By basic inversions, we may assume that the normal orientations $n_1$ and $n_2$  of $e_1$ and $e_2$ point from $f$ to the opposite regions.  
Let $x_1$ and $x_2$ be the labels of them, respectively.  
The {\em sign} of 
$v$ with respect to the region $f$ is $+1$ (or $-1$) if the pair of normal orientations $(n_2, n_1)$ does (or does not) match the orientation of ${\mathbb R}^2$.  
The {\em weight} of $v$ is defined to be $\epsilon (y, x_1, x_2)$ where $\epsilon $ is the sign of $v$. See Figure~\ref{fgweightsa}.

\begin{figure}[H]
\begin{center}
\begin{picture}(240,50)(0,0)

\put(30,0){\line(0,1){50}} 
\put(5,25){\line(1,0){20}}
\put(35,25){\line(1,0){20}} 
\put(17,15){\vector(0,1){20}} 
\put(20,10){\vector(1,0){20}} 
\put(8,5){$[y]$}
\put(12,40){$x_1$} 
\put(45,8){$x_2$} 
\put(65,18){\framebox(50,14){$(y, x_1, x_2)$} }

\put(140,25){\line(1,0){20}} 
\put(170,25){\line(1,0){20}}
\put(165,0){\line(0,1){50}} 
\put(150,35){\vector(0,-1){20}} 
\put(157,40){\vector(1,0){20}} 
\put(145,5){$x_1$}
\put(145,40){$[y]$} 
\put(185,42){$x_2$} 
\put(205,18){\framebox(50,14){$ -(y, x_1, x_2)$} }

\end{picture}
\end{center}
\caption{Weights}
\label{fgweightsa}
\end{figure}

\begin{lemma}
As an element of $C_2^{{\rm R}, \rho}(X)_Y = C_2(X)_Y / D_2^{\rho}(X)_Y$, the weight of $v$ does not depend on the specified region. 
\end{lemma}

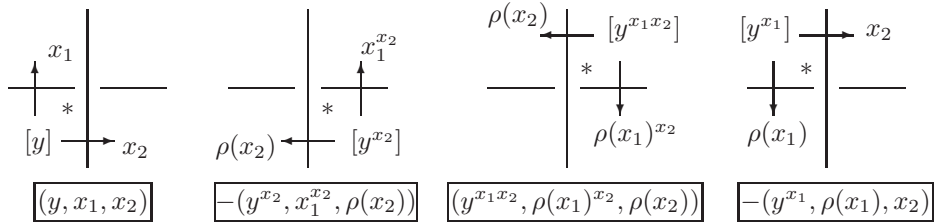
\begin{figure}[H]
\begin{center}
\begin{minipage}{80pt}
\begin{picture}(60,80)
\put(0,50){\line(1,0){25}}
\put(35,50){\line(1,0){25}}
\put(30,20){\line(0,1){60}} 
\put(20,30){\vector(1,0){20}} 
\put(10,40){\vector(0,1){20}}
\put(43,25){$x_2$} 
\put(15,62){$x_1$}
\put(6,27){$[y]$}
\put(20,40){$\ast $}
\put(10,0){\framebox(44,14){$(y,x_1,x_2)$}}
\end{picture}
\end{minipage}
\begin{minipage}{80pt}
\begin{picture}(60,80)
\put(0,50){\line(1,0){25}}
\put(35,50){\line(1,0){25}}
\put(30,20){\line(0,1){60}} 
\put(40,30){\vector(-1,0){20}} 
\put(50,40){\vector(0,1){20}}
\put(-5,25){$\rho(x_2)$} 
\put(50,63){$x_1^{x_2}$}
\put(46,27){$[y^{x_2}]$}
\put(35,40){$\ast $}
\put(-5,0){\framebox(76,14){$-(y^{x_2},x_1^{x_2},\rho(x_2))$}}
\end{picture}
\end{minipage}
\hspace{4mm}
\begin{minipage}{80pt}
\begin{picture}(60,80)
\put(0,50){\line(1,0){25}}
\put(35,50){\line(1,0){25}}
\put(30,20){\line(0,1){60}} 
\put(40,70){\vector(-1,0){20}} 
\put(50,60){\vector(0,-1){20}}
\put(40,30){$\rho(x_1)^{x_2}$} 
\put(0,75){$\rho(x_2)$}
\put(45,70){$[y^{x_1 x_2}]$}
\put(35,55){$\ast $}
\put(-15,0){\framebox(96,14){$(y^{x_1 x_2},\rho(x_1)^{x_2},\rho(x_2))$}}
\end{picture}
\end{minipage}
\hspace{4mm}
\begin{minipage}{80pt}
\begin{picture}(60,80)
\put(0,50){\line(1,0){25}}
\put(35,50){\line(1,0){25}}
\put(30,20){\line(0,1){60}} 
\put(20,70){\vector(1,0){20}} 
\put(10,60){\vector(0,-1){20}}
\put(0,30){$\rho(x_1)$} 
\put(45,70){$x_2$}
\put(-3,70){$[y^{x_1}]$}
\put(20,55){$\ast $}
\put(-3,0){\framebox(73,14){$-(y^{x_1},\rho(x_1),x_2)$}}
\end{picture}
\end{minipage}
\end{center}
\caption{Weights for the same crossing}
\label{fgweightsb}
\end{figure}

\begin{proof}
When we change the specified region as in Figure~\ref{fgweightsb}, where the specified region is denoted by $\ast $,
the weight $(y,x_1,x_2)$ changes to $-(y^{x_2},x_1^{x_2},\rho(x_2))$, $(y^{x_1 x_2},\rho(x_1)^{x_2},\rho(x_2))$ or $-(y^{x_1},\rho(x_1),x_2)$.
They are the same element of $C_2(X)_Y/D_2^{\rho}(X)_Y$.
\end{proof}

For a diagram $D$ and an $(X,\rho)_Y$-coloring $C$, we define a chain $c_{D, C}$ by 
\[
 c_{D, C} = \sum_{v} \epsilon (y, x_1, x_2) \quad \in C_2^{{\rm R}, \rho}(X)_Y  ~ {\rm or}  ~ C_2^{{\rm Q},\rho}(X)_Y, 
\]
where $v$ runs over all crossings of $D$ and $\epsilon (y, x_1, x_2)$ is the weight of $v$. 

Let $(D,C)$ be a pair of a diagram and an $(X,\rho)_Y$-coloring.
If $D$ changes into $D'$ by a single Reidemeister move in a disk support $E$,
then there is a unique $(X,\rho)_Y$-coloring, say $C'$, of $D'$ which is identical with $C$ outside $E$.   In this situation, we say that $(D',C')$ is obtained from $(D,C)$ by a Reidemeister move (with support $E$).

We say that $(D, C)$ and $(D', C')$ are {\it Reidemeister move equivalent} if they are related by a finite sequence of Reidemeister moves.  The equivalence class of $(D, C)$ is called an {\em $(X, \rho)_Y$-colored link}.

\begin{theorem}
\label{theorem:2-cycle} 
 The homology class $[c_{D,C}] \in H_2^{{\rm Q}, \rho}(X)_Y$ is an invariant of an $(X, \rho)_Y$-colored link.  
\end{theorem}
 
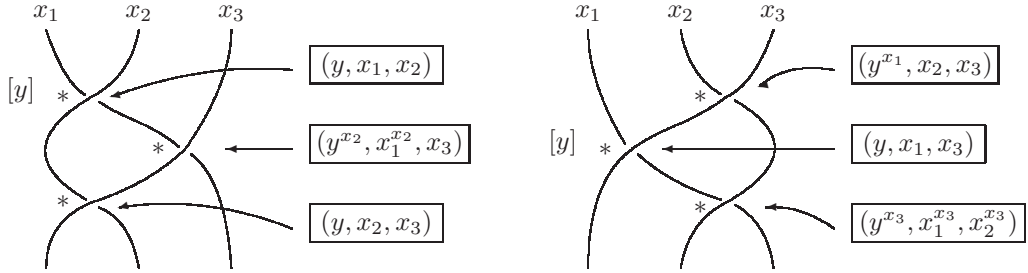
\begin{figure}[H]
\begin{center}
\begin{minipage}{150pt}
\begin{picture}(230,110)
\qbezier(0,90)(7,70)(14.5,65.8)
\qbezier(20,63)(24,60)(35.5,55)
\qbezier(50,46)(46,50)(35.5,55)
\qbezier(55,43)(68,33)(70,0)
   
\qbezier(35,90)(32,73)(17.5,65)
\qbezier(15,26)(-17,45)(17.5,65)
\qbezier(20,23)(32,18)(35,0)

\qbezier(70,90)(70,70)(52.5,45)
\qbezier(17.5,25)(39,32)(52.5,45)
\qbezier(17.5,25)(0,18)(0,0)

\put(100,70){\framebox(50,14){$(y,x_1,x_2)$}}
\put(100,40){\framebox(60,14){$(y^{x_2},x_1^{x_2},x_3)$}}
\put(100,10){\framebox(50,14){$(y,x_2,x_3)$}}

\qbezier(93,75)(60,78)(25,65)
\put(25,65){\vector(-4,-1){}}
\put(93,45){\vector(-1,0){25}}
\qbezier(93,15)(60,28)(28,23)
\put(28,23){\vector(-4,-1){}}

\put(-5,95){$x_1$}
\put(30,95){$x_2$}
\put(65,95){$x_3$}
\put(-15,65){$[y]$}

\put(4,60){*}
\put(4,20){*}
\put(40,40){*}
\end{picture}
\end{minipage}
\hspace{17mm}
\begin{minipage}{150pt}
\begin{picture}(260,110)

\qbezier(70,90)(63,70)(52.5,65)
\qbezier(52.5,65)(46,60)(35.5,55)
\qbezier(16.5,45)(23,50)(35.5,55)
\qbezier(16.5,45)(0,33)(0,0)
   
\qbezier(35,90)(38,73)(50,66)
\qbezier(52.5,25)(87,45)(55.5,63.5)
\qbezier(52.5,25)(38,18)(35,0)

\qbezier(0,90)(0,70)(14,47)
\qbezier(50,26)(31,32)(19,43)
\qbezier(55.5,24)(68,17)(70,0)

\put(100,70){\framebox(55,14){$(y^{x_1},x_2,x_3)$}}
\put(100,40){\framebox(50,14){$(y,x_1,x_3)$}}
\put(100,10){\framebox(65,14){$(y^{x_3},x_1^{x_3},x_2^{x_3})$}}
 
\qbezier(93,75)(75,78)(65,69)
\put(65,69){\vector(-4,-1){}}
\put(93,45){\vector(-1,0){65}}
\qbezier(93,15)(80,24)(68,23)
\put(68,23){\vector(-1,0){}}

\put(-5,95){$x_1$}
\put(30,95){$x_2$}
\put(65,95){$x_3$}
\put(-15,45){$[y]$}

\put(4,40){*}
\put(40,18){*}
\put(40,60){*}

\end{picture}
\end{minipage}
\end{center}
\caption{A Reidemeister move of type III}
\label{Reidemeister move}
\end{figure}

\begin{proof} 
First we show that the chain $c_{D,C}$ is a $2$-cycle of $C_\ast^{{\rm R}, \rho}(X)_Y$ and $C_\ast^{{\rm Q}, \rho}(X)_Y$.   
Let us fix a checkerboard coloring of the complementary regions of $D$. 
Using basic inversions, we assume that normal orientations of semi-arcs of $D$ point from the black regions to the white.
Then we give an orientation to every semi-arc of $D$ such that 
the orientation vector followed by the normal orientation vector matches the orientation of ${\mathbb R}^2$.
Assign a pair $(y,x)\in Y\times X$ to each semi-arc, where $y$ is the color of the black region and $x$ is the color of the semi-arc.
See Figure~\ref{fglabelsemiarc}.

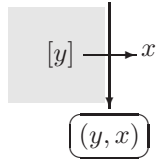
\begin{figure}[H]
\begin{center}
\begin{picture}(50,40)(0,0)

\put(25,40){\vector(0,-1){40}}  
\put(-13,20){\setlength{\fboxrule}{15pt}\fcolorbox[gray]{.9}{.9}{}} 
\put(15,20){\vector(1,0){20}} 
\put(1,18){$[y]$}
\put(37,20){$x$} 
\put(10,-13){\ovalbox{$(y,x)$}}
\end{picture}
\end{center}
\caption{The label of a semi-arc}
\label{fglabelsemiarc}
\end{figure}

If the semi-arc is not a simple loop, then we give the initial point the weight $-(y,x)$ and the terminal point the weight $(y,x)$ with respect to the orientation of the semi-arc.
For each crossing $v$ of $D$, we choose its specified region from the black regions.
See Figure~\ref{fgweightsc}.
If Figure~\ref{fgweightsc}(1) occurs at $v$, the weights of $v$ is $(y,x_1,x_2)$, and
\begin{align*}
\partial _2(y,x_1,x_2)&= -(y,x_2)+(y^{x_1},x_2)+(y,x_1)-(y^{x_2},x_1^{x_2})\\
                      &\equiv  -(y,x_2) -(y^{x_1 x_2},\rho (x_2))+(y,x_1)+ (y^{x_1 x_2},\rho(x_1)^{x_2}) 
\end{align*}
mod $D_1^{\rho}(X)_Y.$
The four terms are exactly the same with the weights assigned the endpoints of semi-arcs. See Figure~\ref{fgweightsd}.  It is similar for Figure~\ref{fgweightsc}(2).  

\begin{figure}[H]
\begin{center}
\begin{picture}(240,50)(0,0)

\put(30,0){\line(0,1){50}} 
\put(5,25){\line(1,0){20}}
\put(35,25){\line(1,0){20}} 
\put(6,10){\setlength{\fboxrule}{8pt}\fcolorbox[gray]{.9}{.9}{}}
\put(35,38){\setlength{\fboxrule}{8pt}\fcolorbox[gray]{.9}{.9}{}}
\put(17,15){\vector(0,1){20}} 
\put(20,10){\vector(1,0){20}} 
\put(20,15){$\ast $}
\put(8,5){$[y]$}
\put(12,40){$x_1$} 
\put(45,8){$x_2$} 
\put(65,18){\framebox(50,14){$(y, x_1, x_2)$} }
\put(25,-15){(1)}

\put(140,25){\line(1,0){20}} 
\put(170,25){\line(1,0){20}}
\put(165,0){\line(0,1){50}} 
\put(140,38){\setlength{\fboxrule}{8pt}\fcolorbox[gray]{.9}{.9}{}}
\put(168,10){\setlength{\fboxrule}{8pt}\fcolorbox[gray]{.9}{.9}{}}
\put(150,35){\vector(0,-1){20}} 
\put(160,40){\vector(1,0){20}} 
\put(144,40){$[y]$}
\put(182,40){$x_2$} 
\put(145,8){$x_1$} 
\put(155,28){$\ast $}
\put(205,18){\framebox(50,14){$ -(y, x_1, x_2)$} }
\put(160,-15){(2)}

\end{picture}
\end{center}
\caption{Weights}
\label{fgweightsc}
\end{figure}

\begin{figure}[H]
\begin{center}
\begin{picture}(130,95)(0,0)
\put(50,45){\circle*{2}}
\put(50,45){\vector(0,1){30}} 
\put(50,45){\vector(0,-1){30}} 
\put(5,45){\vector(1,0){40}}
\put(95,45){\vector(-1,0){40}}
\put(5,30){\setlength{\fboxrule}{10pt}\fcolorbox[gray]{.9}{.9}{\ \ \ \ \ }}
\put(52,60){\setlength{\fboxrule}{10pt}\fcolorbox[gray]{.9}{.9}{\ \ \ \ \ }} 
\put(30,35){\vector(0,1){20}} 
\put(40,25){\vector(1,0){20}} 
\put(60,65){\vector(-1,0){20}} 
\put(70,55){\vector(0,-1){20}} 
\put(26,25){$[y]$}
\put(15,50){$x_1$} 
\put(55,15){$x_2$} 
\put(65,62){$[y^{x_1 x_2}]$}
\put(25,70){$\rho (x_2)$}
\put(76,33){$\rho (x_1)^{x_2}$}
\put(110,39){\ovalbox{$+(y^{x_1 x_2},\rho(x_1)^{x_2})$}}
\put(-50,39){\ovalbox{$+(y,x_1)$}}
\put(33,-5){\ovalbox{$-(y,x_2)$}}
\put(22,88){\ovalbox{$-(y^{x_1 x_2},\rho (x_2))$}}
\end{picture}
\end{center}
\caption{Weights}
\label{fgweightsd}
\end{figure}

Since the weights of endpoints of semi-arcs are canceled, we see that $\partial _2(c_{D,C})\equiv 0 \mod{D_1^{\rho}(X)_Y}$.
Therefore the chain $c_{D,C}$ is a $2$-cycle of $C_\ast^{{\rm R}, \rho}(X)_Y$ and $C_\ast^{{\rm Q},\rho}(X)_Y$. 

It is sufficient to prove that 
if $(D, C)$ changes into $(D', C')$ by a single Reidemeister move 
then 
$[c_{D, C}]= [c_{D', C'}]$ in  $H_2^{{\rm Q}, \rho}(X)_Y$. 
If the Reidemeister move is of type II,
then $c_{D,C}$ equals $c_{D',C'}$ in the chain group $C_2^{{\rm R},\rho}(X)_Y$.
If the Reidemeister move is of type III (Figure~\ref{Reidemeister move}),
then by choosing specified regions as in Figure~\ref{Reidemeister move}, it is seen that the difference of $c_{D,C}$ and $c_{D',C'}$ 
is $\pm \partial _3(y,x_1,x_2,x_3)$ for some $(y,x_1,x_2,x_3)\in Y\times X^3$, cf. \cite{CJKLS03, FRS95, RS1, RS}.
In either case, $[c_{D,C}]=[c_{D',C'}]$ in $H_2^{{\rm R},\rho}(X)_Y$ and in $H_2^{{\rm Q},\rho}(X)_Y$.
If the Reidemeister move is of type I, then the difference of $c_{D,C}$ and $c_{D',C'}$ is $\pm (y,x,x)$ for some $(y,x)\in Y\times X$,
which vanishes in $C_2^{{\rm Q},\rho}(X)_Y$. \end{proof}

Let 
$$
{\cal H}(D) = \{ [c_{D, C}] \in H_2^{{\rm Q}, \rho}(X)_Y \, 
| \, C: \mbox{ $(X, \rho)_Y$-colorings of $D$} \} 
$$
as a multi-set.   (A {\it multi-set} means a set with (possible) repeats.)

For  a $2$-cocycle $\theta$ of the cochain complex $C^*_{{\rm Q}, \rho}(X, A)_Y$ with a coefficient group $A$, let 
$$
\Phi _\theta(D) = \{ \theta(c_{D, C}) \in A \, 
| \, C: \mbox{ $(X, \rho)_Y$-colorings of $D$} \} 
$$
as a multi-set.  

\begin{corollary}\label{theorem:invariant}
The multi-sets ${\cal H}(D)$ and $\Phi _\theta(D)$ are invariants of the link type of $D$.  
\end{corollary} 

By  definition of $C_{{\rm Q},\rho}^*(X, A)_Y$, the linear extension $\theta : \mathbb{Z}(Y\times X^2)\to A$ 
of a map $\theta : Y\times X^2 \to A$ is a $2$-cocycle 
of $C_{{\rm Q},\rho}^*(X, A)_Y$ 
if and only if the following conditions {\rm (1)--(3)} are satisfied. 
\begin{description}
\item[(1)] For any $(y,x_1,x_2,x_3)\in Y\times X^3$, 
\begin{multline*}
- \theta(y,x_2,x_3) + \theta(y^{x_1},x_2,x_3) + \theta(y,x_1,x_3)\\
- \theta(y^{x_2},x_1^{x_2},x_3) -  \theta(y,x_1,x_2) + \theta(y^{x_3},x_1^{x_3},x_2^{x_3}) =0,
\end{multline*}
\item[(2)] For any $(y,x)\in Y\times X$, $\theta(y,x,x)=0$, and
\item[(3)] For any $(y,x_1,x_2)\in Y\times X^2$,
\[
\begin{array}{l}
\theta(y,x_1,x_2 )+ \theta(y^{x_1}, \rho(x_1),x_2)=0\quad{\rm and} \quad 
\theta(y,x_1,x_2) + \theta(y^{x_2},x_1^{x_2}, \rho(x_2)) =0.
\end{array}
\]
\end{description}
We call these conditions the {\em symmetric quandle $2$-cocycle conditions}.

\begin{example}\label{example:linking}{\rm 
Let $X$ be 
the order $4$ trivial quandle $T_4=  \{e_1,e'_1, e_2,e'_2\}$,  and  
$\rho: X \to X$ a good involution with $\rho(e_i) = e'_i$ $(i =1,2)$ (see Proposition~\ref{proposition:trivial}).  Let $Y=\{e\}$, which is an $(X, \rho)$-set. 
Define a map $\theta : Y\times X^2 \to {\mathbb Z}$ by 
$$  \theta = \chi_{(e, e_1, e_2)} + \chi_{(e, e'_1, e'_2)}  
- \chi_{(e, e'_1, e_2)} - \chi_{(e, e_1, e'_2)},  
$$
where $\chi_{(e, a, b)}$ is defined by 
$\chi_{(e,a,b)} (e,x,y) = 1$ if $(e, x ,y)=(e, a, b)$ and $\chi_{(e,a,b)} (e,x,y) = 0$ otherwise. 
Then $\theta$ satisfies the symmetric quandle $2$-cocycle conditions, and 
the linear extension $\theta : \mathbb{Z}(Y\times X^2)\to {\mathbb Z}$ is a $2$-cocycle. 

Let $D$ be a diagram of an unoriented $(2,m)$-torus link with $2m$ crossings.   
There are 16 $(X, \rho)_Y$-colorings, and they represent cycles 
$$ m (e, a, b ) + m (e, b, a ) $$ 
for $a, b \in X$.  
When we evaluate these cycles by $\theta$, we have 
$$ \theta( m (e, a, b ) + m (e, b, a ) ) = 
\left\{ \begin{array}{ll}
m & \mbox{if } (a,b) = (e_i, e_j), (e'_i, e'_j)   \mbox{ for } i \neq j, \\
-m & \mbox{if } (a,b) = (e'_i, e_j), (e_i, e'_j)  \mbox{ for } i \neq j, \\
0 & \mbox{otherwise }. 
\end{array} \right. 
$$
Thus $\Phi _f(D) =\{ m, m, m, m, -m, -m, -m, -m, 0, \cdots, 0 (\mbox{$8$ times}) \}$.  
}\end{example}

\begin{proposition}
Let $(X, \rho)$, $Y$, and $f$ be as in Example~$\ref{example:linking}$.  
Let $D$ be a diagram of a $2$-component unoriented link $L$. 
Then $\Phi _\theta(D) =\{ m, m, m, m, -m, -m, -m, -m, 0, \cdots, 0 (\mbox{$8$ times})  \}$, 
where $m$ is the linking number of $D$ when we give an orientation to $D$. 
\end{proposition}

\begin{proof}
Give an orientation to $D= D_1 \cup D_2$, and denote this oriented diagram by $D^{+}$.  The semi-arcs of $D$ can be assigned a normal orientation such that the orientation vector of $D^{+}$ followed by the normal orientation vector matches the orientation of ${\mathbb R}^2$.   
For each $(X, \rho)_Y$-coloring, we take a representative such that each semi-arc has a normal vector  this way.  Then each semi-arc is labeled by an element of $X$.  The complementary regions are always colored by $e \in Y$. 
  (This coloring is the $X_Y$-coloring of an oriented link diagram in the sense of 
\cite{CJKLS03, RS1}, etc..)  Since $X$ is a trivial quandle, the labels of semi-arcs of the same knot component $D_k$ $(k=1,2)$ are the same.  Let $x_k$ be the label of $D_k$.   There are 16 admissible colorings according to $(x_1, x_2) \in X \times X$.  
The linking number ${\rm Lk}(D^{+})$ of $D^{+}$ is the sum of signs of crossings where $D_i$ is lower and $D_j$ is upper for $i \neq j$. 
When $(x_1, x_2)= (e_i, e_j)$ or $(e'_i, e'_j)$ with $i \neq j$, then 
the value $\theta(e, a, b)$ at a crossing $v$ is the sign of $v$ if $D_i$ is lower and $D_j$ is upper, or $0$ otherwise.  Thus for these colorings $C$, $\theta(c_{D,C}) = {\rm Lk}(D^{+})$. 
When $(x_1, x_2)= (e'_i, e_j)$ or $(e_i, e'_j)$ with $i \neq j$, then 
the value $\theta(e, a, b)$ at a crossing $v$ is the negative of the sign of $v$ if $D_i$ is lower and $D_j$ is upper, or $0$ otherwise.  Thus for these colorings $C$, $\theta(c_{D,C}) = -{\rm Lk}(D^{+})$. 
When $(x_1, x_2)$ is not $  (\pm e_i, \pm e_j)$ with $i \neq j$, then the value $\theta(e, a, b)$ at a crossing $v$ is $0$.  Thus for such a coloring $C$, $\theta(c_{D,C}) =0$.   \end{proof}

The invariants ${\cal H}(D)$ and $\Phi _\theta(D)$ are analogies of invariants of oriented links in the sense of 
\cite{CJKLS03, CKS01, RS1}, which we denote by ${\cal H}^{\rm ori}$ and $\Phi _{\theta }^{\rm ori}$. We recall them. 
Let $D^+$ be an oriented diagram such that $D$ is obtained from $D^+$ by forgetting the orientation. 
Using the orientation, we assign 
the semi-arcs of $D^+$ normal orientations such that the orientation vector followed by the normal orientation vector matches the orientation of ${\mathbb R}^2$.  An {\it $X_Y$-coloring} of $D^+$  is assignment of an element of $X$ to each semi-arc and an element of $Y$ to each complementary region of $D^+$ satisfying  the conditions illustrated in 
the left of Figure~\ref{fgadmcolora} and in Figure~\ref{fgadmcolorb}.   The weight of each crossing is defined as in Figure~\ref{fgweightsa}.  Let $c_{D^+, C}$ be the sum of the weights of all crossings of $D$.  It is a $2$-cycle of $C_\ast^{\rm Q}(X)_Y$, and  
the homology class  $[c_{D^+, C}] \in H_2^{{\rm Q}}(X)_Y$ is a colored link invariant \cite{CJKLS03, CKS01, RS1}.  
Let 
$
{\cal H}^{\rm ori}(D^+) = \{ [c_{D^+, C}] \in H_2^{{\rm Q}}(X)_Y \, 
| \, C: \mbox{$X_Y$-colorings of $D^+$} \} 
$
and 
$
\Phi _\theta^{\rm ori}(D^+) = \{ \theta(c_{D^+, C}) \in A \, 
| \, C: \mbox{$X_Y$-colorings of $D^+$} \} 
$
as multi-sets, where $\theta$ is a $2$-cocycle of $C^*_{{\rm Q}}(X, A)_Y$.  Then 
${\cal H}^{\rm ori}(D^+)$ and  $\Phi _\theta^{\rm ori}(D^+)$ are oriented link invariants. 

Note that if $\theta : \mathbb{Z}(Y\times X^2)\to A$ is a $2$-cocycle of $C_{{\rm Q},\rho}^*(X, A)_Y$, then $\theta$ is a $2$-cocycle of $C_{{\rm Q}}^*(X, A)_Y$. 

\begin{theorem}\label{thm:equality1dim}
Let $(X, \rho)$ be a symmetric quandle, $Y$ an $(X, \rho)$-set, $A$ an abelian group, and $\theta$ a $2$-cocycle of $C_{{\rm Q},\rho}^*(X, A)_Y$.  Let $D$ be a diagram of an unoriented link in ${\mathbb R}^3$.  Let $D^+$ be an arbitrarily oriented diagram of $D$.  
Then $$ 
\Phi _\theta(D)= \Phi^{\rm ori}_\theta(D^+). $$
\end{theorem}

\begin{proof}
Assign a normal orientation to each semi-arc of $D$ (and of $D^+$) such that the orientation vector determined from $D^+$ followed by the normal orientation vector matches the orientation of ${\mathbb R}^2$.  
For each  $(X, \rho)_Y$-coloring $C$ of $D$, we take a representative of $C$ such that the normal orientation of each semi-arc is the same with the normal orientation.  Then $C$ determines an assignment of an 
element of $X$ to each semi-arc.  
Let $C^+$ be this assignment together with the same assignment of an element of $Y$ to each complementary region as $C$.   Then $C^+$ is an $X_Y$-coloring for the oriented diagram $D^+$ in the sense of \cite{CJKLS03, RS1}.  Thus there is a bijection between the set of $(X, \rho)_Y$-colorings of $D$ and the set of  
$X_Y$-colorings of $D^+$.  At each crossing,   the weights  for $(D,C)$ and $(D^+,C^+)$ evaluated by $\theta$ are the same, and hence we have 
$$ \theta(c_{D, C}) = \theta(c_{D^+, C^+}) $$
and  $ 
\Phi_\theta(D^+) =  \Phi^{\rm ori} _\theta(D)$. 
\end{proof}

\begin{corollary}\label{corollary:equality1dim}
Let $(X, \rho)$ be a symmetric quandle, $Y$ an $(X, \rho)$-set, $A$ an abelian group, and $\theta$ a $2$-cocycle of $C_{{\rm Q}}^*(X, A)_Y$.  The invariant $\Phi^{\rm ori}_\theta$ does not depend on the orientation of a link if $\theta$ is cohomologous in $C_\ast^{\rm Q}(X, A)_Y$ 
to a $2$-cocycle satisfying the symmetric quandle $2$-cocycle conditions. 
\end{corollary}

\section{A quick calculation of $\Phi_\theta$}\label{Link_invariants_2}
Throughout this section, we assume that $X$ is a kei and $\rho$ is the identity map of $X$ (cf. Proposition~\ref{kei-goodinvolution}).
Under this assumption, the definition of the invariants $\Phi_\theta$ can be much simplified.

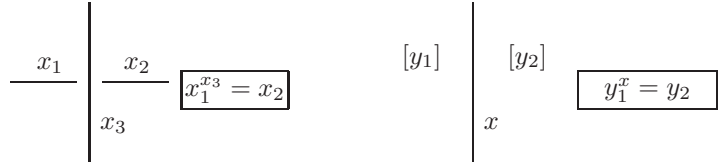
\begin{figure}[H]
\begin{center}
\begin{minipage}{110pt}
\begin{picture}(110,60)
\put(0,30){\line(1,0){25}}
\put(35,30){\line(1,0){25}}
\put(30,0){\line(0,1){60}} 
\put(34,12){$x_3$} 
\put(43,35){$x_2$} 
\put(10,35){$x_1$}
\put(65,20){\framebox(40,14){$x_1^{x_3} = x_2$}}
\end{picture}
\end{minipage}
\hspace{10mm}
\begin{minipage}{110pt}
\begin{picture}(110,60)
\put(30,0){\line(0,1){60}} 
\put(34,12){$x$} 
\put(43,37){$[y_2]$} 
\put(3,37){$[y_1]$}
\put(70,20){\framebox(52,14){$y_1^{x} =y_2$}}
\end{picture}
\end{minipage}
\end{center}
\caption{Coloring conditions}
\label{fgadmcolorc}
\end{figure}

Let $D$ be a diagram in ${\mathbb R}^2$ of an unoriented link in ${\mathbb R}^3$. 
An {\em $X$-coloring} of $D$ is an assignment of an element of $X$ to each arc  of $D$ such that 
for each crossing of $D$, $x_1^{x_3}=x_2$ holds where $x_1$ and $x_2$ are elements of $X$ assigned under-arcs
and $x_3$ is the element assigned the over-arc as in Figure~\ref{fgadmcolorc}.  

Let $Y$ be an $X$-set such that $y^x = y^{x^{-1}}$ for any $(x,y) \in X \times Y$. An {\em $X_Y$-coloring} of $D$ is an 
$X$-coloring of $D$ with an assignment of an element of $Y$ to each complementary region of $D$ such that,  for each arc of $D$,
$y_1^x=y_2$ holds where $y_1$ and $y_2$ are elements of $Y$ assigned the regions separated by the arc, and $x$ is the element assigned the arc as in Figure~\ref{fgadmcolorc}.

Since $X$ is a kei and $\rho={\rm id}_X$, we see the following.
\begin{description}
\item[(1)] An $X$-set  $Y$ with $y^x = y^{x^{-1}}$ for any $(x,y) \in X \times Y$ is an $(X,\rho)$-set, and vice versa.
\item[(2)] Forgetting normal orientations, we have a bijection from the set of $(X,\rho)$-colorings of $D$ to the set of $X$-colorings of $D$.
\item[(3)] Forgetting normal orientations, we have a bijection from the set of $(X,\rho)_Y$-colorings of $D$ to the set of $X_Y$-colorings of $D$.
\end{description}

The following is well-known, which is a special case of Proposition~\ref{proposition:coloring}.

\begin{proposition}{\rm (cf. \cite{CJKLS03,FR92,Joyce, Matveev,RS}) }
Let $X$ be a kei, and $Y$ an $X$-set  such that $y^x = y^{x^{-1}}$ for any $(x,y) \in X \times Y$. 
If two link diagrams represent the same unoriented link type, then there is a bijection between the sets of $X$-colorings of the diagrams, and there is a bijection between the sets of $X_Y$-colorings.
\end{proposition}

Let $D$ be an unoriented link diagram.
Fix an $X_Y$-coloring of $D$, say $C$.
For a crossing $v$ of $D$, there are four complementary regions of $D$ around $v$.
Choose one of them, and let $y$ be the element of $Y$ assigned the region.
Let $x_1$ and $x_2$ be elements of $X$ assigned the under and over-arcs facing the region as in Figure~\ref{fgweightse}.  
The {\em weight} of $v$ is $\epsilon (y, x_1, x_2)$ where $\epsilon $ is the sign of $v$.
(Recall that the sign of a crossing of an unoriented link diagram is defined when a region around the crossing is specified.)  

\begin{figure}[H]
\begin{center}
\begin{picture}(240,50)(0,0)

\put(30,0){\line(0,1){50}} 
\put(5,25){\line(1,0){20}}
\put(35,25){\line(1,0){20}} 
\put(10,10){$[y]$}
\put(11,30){$x_1$} 
\put(33,8){$x_2$} 
\put(65,18){\framebox(50,14){$(y, x_1, x_2)$} }

\put(140,25){\line(1,0){20}} 
\put(170,25){\line(1,0){20}}
\put(165,0){\line(0,1){50}} 
\put(145,16){$x_1$}
\put(145,36){$[y]$} 
\put(168,38){$x_2$} 
\put(205,18){\framebox(50,14){$ -(y, x_1, x_2)$} }

\end{picture}
\end{center}
\caption{Weights}
\label{fgweightse}
\end{figure}

Now we obtain the chain $c_{D, C}$ as the sum of the weights of all crossing of $D$.   
The following is a special case of Theorem~\ref{theorem:2-cycle} where $X$ is a kei and $\rho$ is the identity map. 

\begin{theorem} 
\label{theorem:2-cyclekeiversion}
The homology class $ [c_{D,C} ] \in H_2^{{\rm Q}, \rho}(X)_Y$ is an invariant of an $X_Y$-colored link. 
\end{theorem}

Let 
${\cal H}(D) = \{ [c_{D, C}] \in H_2^{{\rm Q}, \rho}(X)_Y \, 
| \, C: \mbox{$X_Y$-colorings of $D$} \}$ and 
$\Phi _\theta(D) = \{ \theta(c_{D, C}) \in A \, 
| \, C: \mbox{$X_Y$-colorings of $D$} \} $ 
 as multi-sets,  where $\theta$ is a 
$2$-cocycle  of the cochain complex $C^*_{{\rm Q}, \rho}(X, A)_Y$.  

\begin{corollary}\label{theorem:invariantkeiversion}
The multi-sets ${\cal H}(D)$ and $\Phi _\theta(D)$ are invariants of the link type of $D$.  
\end{corollary}

In order to obtain an invariant of $D$ as in Corollary~\ref{theorem:invariantkeiversion}, one may restrict the range of colorings $C$ suitably in the definition of ${\cal H}(D)$ or $\Phi _\theta(D)$;
for example, we have another invariant by assuming that $C$ runs over all colorings of $D$ such that the unbounded region is colored by a specified element of $Y$ (cf. Example~\ref{example:trefoil}). 

\begin{example}{\rm 
\label{example:trefoil}
Let $X$ be the dihedral quandle $R_3= \mathbb{Z}/3\mathbb{Z}$ of order $3$, and $\rho$ the identity map of $X$. 
Let $Y=X$ on which $X$ acts by the quandle operation. 
A map 
$$
\theta : X \times X^2 \to \mathbb{Z}/3\mathbb{Z} \, ; \/   \theta(x,y,z)=(x-y)(y-z)^2z
$$ 
satisfies the symmetric quandle $2$-cocycle conditions  and the linear extension 
$\theta : \mathbb{Z}(X\times X^2)\to \mathbb{Z}/3\mathbb{Z}$ 
is a $2$-cocycle of the cochain complex $C_{{\rm Q},\rho}^*(X, \mathbb{Z}/3\mathbb{Z})_X$.
(It is a modified version of Mochizuki's cocycle \cite{M}.) 

Thus, for an unoriented link diagram $D$ in $\mathbb{R}^2$, a multi-set 
\[
\Phi _\theta(D)_{\infty =0}=\left\{
\begin{array}{l|l}
\theta(c_{D,C})\in \mathbb{Z}/3\mathbb{Z}& C:\mbox{$X_X$-colorings of {\it D} such that}\\
&\mbox{the unbounded region is colored by $0$} 
\end{array}
\right\}
\]
is an unoriented link invariant. 
For example, let $D$ be the diagram of a trefoil illustrated in Figure~\ref{example1}.
Any coloring $C$ of $D$ such that the unbounded region is $0$ is given as in Figure~\ref{example1},
where $a,b \in X =\mathbb{Z}/3\mathbb{Z}$,
and the $2$-chain derived from the colored diagram $(D,C)$ is 
\[
c_{D,C}=(0,a,b)+(0,b,-a-b)+(0,-a-b,a).
\]

\begin{figure}
\begin{center}
\begin{picture}(300,150)(-20,-35)
 \qbezier(19,50)(17,13)(50,-3)
 \qbezier(19,54)(50,70)(78,55)
 \qbezier(55,-1)(80,15)(81,54)

 \qbezier(-5,5)(6,-15)(47,-4)
 \qbezier(-5,5)(-17,34)(19,54)

 \qbezier(105,5)(94,-15)(50,-3)
 \qbezier(105,5)(117,34)(84,53)

 \qbezier(50,92)(80,90)(81,54)
 \qbezier(50,92)(20,90)(19,57)

 \put(180,-30){\line(1, 1){30}}
 \put(210,-30){\line(-1,1){13}}
 \put(193,-13){\line(-1,1){13}}

 \put(180,30){\line(1,1){30}}
 \put(210,30){\line(-1,1){13}}
 \put(193,47){\line(-1,1){13}}

 \put(180,90){\line(1,1){30}}
 \put(210,90){\line(-1,1){13}}
 \put(193,107){\line(-1,1){13}}

\put(145,-35){$-a-b$}
\put(166,-17){$[0]$}
\put(169,0){$b$}

\put(169,25){$a$}
\put(166,42){$[0]$}
\put(145,60){$-a-b$}

\put(169,85){$b$}
\put(166,102){$[0]$}
\put(169,120){$a$}

\put(230,-20){\framebox(58,14){ $(0,b,-a-b)$} }
\put(230,40){\framebox(58,14){$(0,-a-b,a)$} }
\put(230,100){\framebox(37,14){$(0,a,b)$}}

\put(33,30){$[a-b]$}
\put(0,7){$[-b]$}
\put(74,7){$[a+b]$}
\put(38,70){$[-a]$}
\put(-20,65){$[0]$}
\put(-15,4){$b$}
\put(50,100){$a$}
\put(105,-5){$-a-b$}

\put(69,-22){\vector(-1,1){15}}
\qbezier(69,-22)(75,-30)(155,-17)

\put(155,38){\vector(-4,1){65}}

\put(4,85){\vector(1,-3){7}}
\qbezier(4,85)(0,159)(155,102)

\end{picture}
\caption{A diagram of a right-handed trefoil}
\label{example1}
\end{center}
\end{figure}
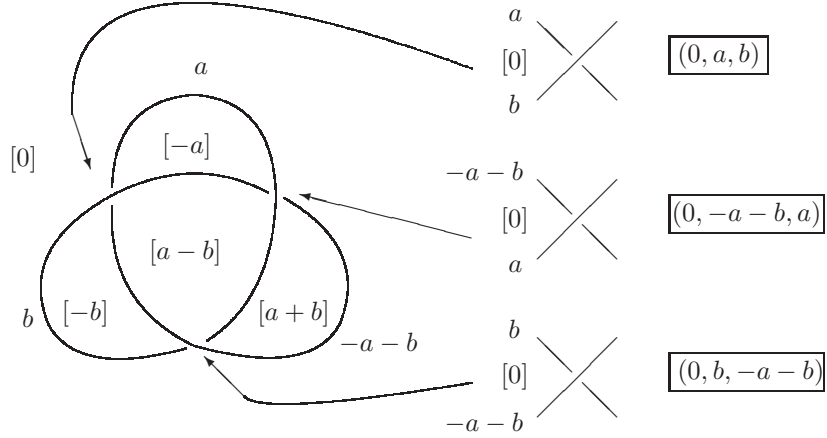

Then we have 
\[
\theta(c_{C,D})  \\
=\begin{cases}
0,&{\rm if}\ a=b,\\
1,&{\rm if}\ a\not =b,
\end{cases}
\]
and hence $\Phi _\theta(D)_{\infty =0}=\{ 0,0,0,1,1,1,1,1,1 \}$,
which is regarded as a multi-set with elements of $\mathbb{Z}/3\mathbb{Z}$.
Let $D^*$ be the mirror image of $D$.
Any coloring $C^*$ of $D^*$ such that the unbounded region is $0$ is given as in Figure~\ref{example2},
and the $2$-chain derived from the colored diagram $(D^*,C^*)$ is 
\[
c_{D^*,C^*}=-(0,a,b)-(0,b,-a-b)-(0,-a-b,a)=-c_{D,C}.
\] 

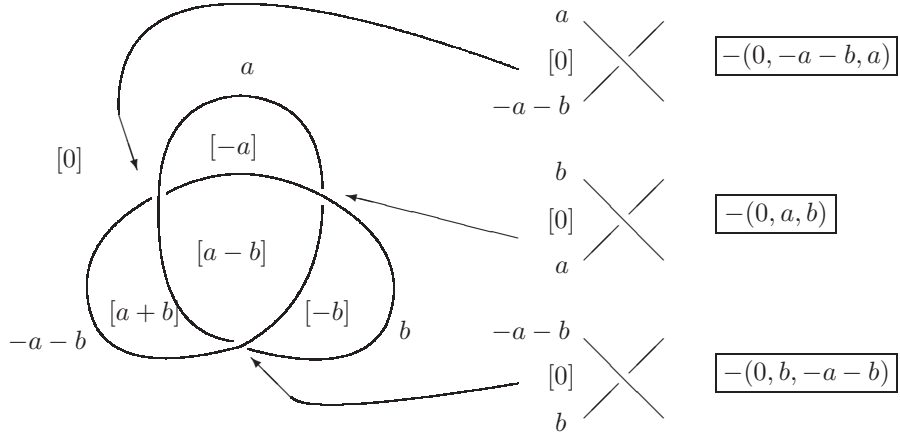
\begin{figure}
\begin{center}
\begin{picture}(300,160)(-10,-10)
 \qbezier(19,74)(17,23)(47,19)
 \qbezier(22,75)(50,90)(81,74)
 \qbezier(50,17)(80,33)(81,70)

 \qbezier(-5,25)(6,5)(50,17)
 \qbezier(-5,25)(-17,54)(16,73)

 \qbezier(105,25)(94,5)(53,16)
 \qbezier(105,25)(117,54)(81,74)

 \qbezier(50,112)(80,110)(81,77)
 \qbezier(50,112)(20,110)(19,74)

 \put(180,20){\line(1, -1){30}}
 \put(180,-10){\line(1,1){13}}
 \put(197,7){\line(1,1){13}}

 \put(180,80){\line(1,-1){30}}
 \put(180,50){\line(1,1){13}}
 \put(197,67){\line(1,1){13}}

 \put(180,140){\line(1,-1){30}}
 \put(180,110){\line(1,1){13}}
 \put(197,127){\line(1,1){13}}

\put(169,-15){$b$}
\put(166,3){$[0]$}
\put(145,20){$-a-b$}

\put(169,45){$a$}
\put(166,62){$[0]$}
\put(169,80){$b$}

\put(145,105){$-a-b$}
\put(166,122){$[0]$}
\put(169,140){$a$}

\put(230,0){\framebox(68,14){$-(0,b,-a-b)$} }
\put(230,60){\framebox(45,14){$-(0,a,b)$} }
\put(230,120){\framebox(68,14){$-(0,-a-b,a)$} }

\put(33,50){$[a-b]$}
\put(0,27){$[a+b]$}
\put(74,27){$[-b]$}
\put(38,90){$[-a]$}
\put(-20,85){$[0]$}
\put(-38,16){$-a-b$}
\put(50,120){$a$}
\put(110,20){$b$}

\put(69,-2){\vector(-1,1){15}}
\qbezier(69,-2)(75,-10)(155,3)

\put(155,58){\vector(-4,1){65}}

\put(4,105){\vector(1,-3){7}}
\qbezier(4,105)(0,179)(155,122)

\end{picture}
\caption{A diagram of a left-handed trefoil}
\label{example2}
\end{center}
\end{figure}

Thus we have $\Phi _\theta(D^*)_{\infty =0}=\{ 0,0,0,-1,-1,-1,-1,-1,-1 \}$.
As a consequence, we see that the trefoil is chiral.
This argument is an unoriented link version of the argument given by Rourke and Sanderson \cite{RS1,RS}.
}
\end{example}

As seen in Example~\ref{example:trefoil}, we can calculate the quandle cocycle invariant 
$\Phi^{\rm ori}_\theta$ by calculating $\Phi_\theta$ 
without considering the sign and the orientation of arcs at each crossing if $X$ is a kei and $\rho$ is the identity map. 
This gives a practical method of calculation of some $\Phi^{\rm ori}_\theta$, cf. \cite{I, SatohA}.


\section{Quandle cocycle invariants of surface-links}
\label{Surface-link invariants}

Let $D$ be a diagram in ${\mathbb R}^3$ of an unoriented surface-link in ${\mathbb R}^4$.  
Divide over-sheets at the double curves and we call the sheets of the result {\em semi-sheets} of $D$. 
Each semi-sheet is a compact orientable surface in ${\mathbb R}^3$ (cf. \cite{Kamada2001}).  

\begin{figure}[H]
\begin{center}
\begin{minipage}{190pt}
\begin{picture}(170,130)

\qbezier(7,47)(7,47)(43,81)
\qbezier(7,47)(7,47)(46,47)
\qbezier(46,47)(46,47)(50,51)
\qbezier(50,81)(43,81)(43,81)

\qbezier(131,81)(131,81)(95,47)
\qbezier(53,47)(53,47)(95,47)
\qbezier(53,47)(53,47)(89,81)
\qbezier(131,81)(89,81)(89,81)

\qbezier(50,8)(50,8)(50,90)
\qbezier(84,122)(50,90)(50,90)
\qbezier(50,8)(50,8)(84,40)
\qbezier(84,47)(84,40)(84,40)
\qbezier(84,122)(84,78)(84,78)

\put(65,85){\vector(1,0){25}}
\put(65,35){\vector(1,0){25}}
\put(87,91){$x_2$}
\put(87,25){$x_1$}
\put(130,50){\framebox(35,14){$x_1 = x_2$} }

\end{picture}
\end{minipage}
\begin{minipage}{190pt}
\begin{picture}(170,130)

\qbezier(7,47)(7,47)(43,81)
\qbezier(7,47)(7,47)(46,47)
\qbezier(46,47)(46,47)(50,51)
\qbezier(50,81)(43,81)(43,81)

\qbezier(131,81)(131,81)(95,47)
\qbezier(53,47)(53,47)(95,47)
\qbezier(53,47)(53,47)(89,81)
\qbezier(131,81)(89,81)(89,81)

\qbezier(50,8)(50,8)(50,90)
\qbezier(84,122)(50,90)(50,90)
\qbezier(50,8)(50,8)(84,40)
\qbezier(84,47)(84,40)(84,40)
\qbezier(84,122)(84,78)(84,78)

\put(50,85){\vector(-1,0){15}}
\multiput(65,85)(-4,0){4}{\line(-2,0){2}}
\put(65,35){\vector(1,0){25}}
\put(26,91){$x_2$}
\put(87,25){$x_1$}
\put(130,50){\framebox(48,14){$x_1 = \rho (x_2)$} }

\end{picture}
\end{minipage}\\
\begin{minipage}{190pt}
\begin{picture}(190,130)

\qbezier(7,47)(7,47)(43,81)
\qbezier(7,47)(7,47)(46,47)
\qbezier(46,47)(46,47)(50,51)
\qbezier(50,81)(43,81)(43,81)

\qbezier(131,81)(131,81)(95,47)
\qbezier(53,47)(53,47)(95,47)
\qbezier(53,47)(53,47)(89,81)
\qbezier(131,81)(89,81)(89,81)

\qbezier(50,8)(50,8)(50,90)
\qbezier(84,122)(50,90)(50,90)
\qbezier(50,8)(50,8)(84,40)
\qbezier(84,47)(84,40)(84,40)
\qbezier(84,122)(84,78)(84,78)

\put(65,35){\vector(1,0){25}}
\put(36,63){\vector(0,1){25}}
\put(100,63){\vector(0,1){25}}
\put(20,88){$x_1$}
\put(107,88){$x_2$}
\put(87,25){$x_3$}
\put(130,50){\framebox(46,14){$x_1^{x_3} = x_2$} }

\end{picture}
\end{minipage}
\begin{minipage}{190pt}
\begin{picture}(190,130)

\qbezier(7,47)(7,47)(43,81)
\qbezier(7,47)(7,47)(46,47)
\qbezier(46,47)(46,47)(50,51)
\qbezier(50,81)(43,81)(43,81)

\qbezier(131,81)(131,81)(95,47)
\qbezier(53,47)(53,47)(95,47)
\qbezier(53,47)(53,47)(89,81)
\qbezier(131,81)(89,81)(89,81)

\qbezier(50,8)(50,8)(50,90)
\qbezier(84,122)(50,90)(50,90)
\qbezier(50,8)(50,8)(84,40)
\qbezier(84,47)(84,40)(84,40)
\qbezier(84,122)(84,78)(84,78)

\put(65,35){\vector(1,0){25}}
\put(36,63){\vector(0,1){25}}
\put(100,52){\vector(0,-1){15}}
\multiput(100,63)(0,-4){3}{\line(0,-2){2}}
\put(20,88){$x_1$}
\put(107,39){$x_2$}
\put(87,25){$x_3$}
\put(130,50){\framebox(53,14){$x_1^{x_3} = \rho(x_2)$} }

\end{picture}
\end{minipage}

\end{center}
\caption{Coloring conditions}
\label{Fig.C}
\end{figure}
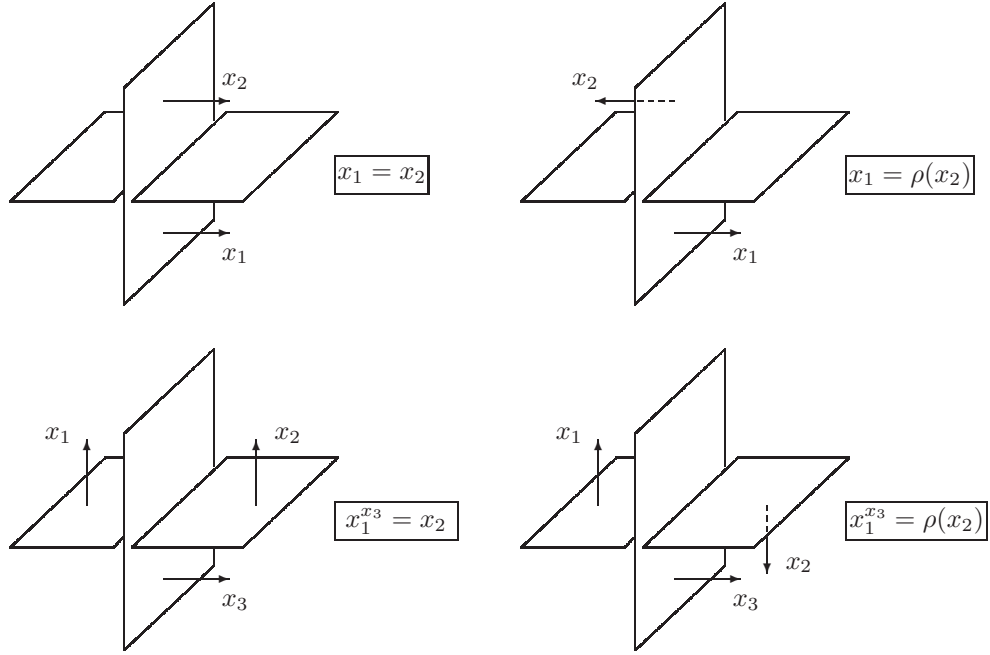

We say that an assignment of a normal orientation and an element of $X$ to each semi-sheet of $D$ 
satisfies the {\em coloring conditions} if it satisfies the following conditions (Figure~\ref{Fig.C}): 
\begin{itemize}
\item 
Suppose that two diagonal semi-sheets coming from an over-sheet of $D$ about a double curve are labeled by $x_1$ and $x_2$. If the normal orientations are coherent then $x_1=x_2$, otherwise $x_1= \rho(x_2)$. 
\item 
Suppose that two diagonal semi-sheets $e_1$ and $e_2$ which are under-sheets about a double curve are labeled by $x_1$ and $x_2$, and suppose that one of the two semi-sheets coming from an over-sheet of $D$, say $e_3$, is labeled by $x_3$.  We assume that the normal orientation of the over semi-arc $e_3$ points from $e_1$ to $e_2$.  If the normal orientations of $e_1$ and $e_2$ are coherent, then $x_1^{x_3}= x_2$, otherwise $x_1^{x_3}= \rho(x_2)$. 
\end{itemize}

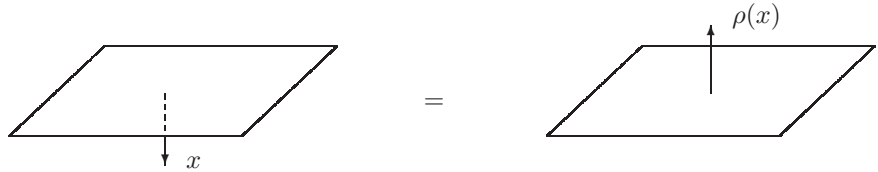
\begin{figure}[H]
\begin{center}
\begin{minipage}{132pt}
\begin{picture}(132,70)

\qbezier(7,22)(7,22)(95,22)
\qbezier(7,22)(7,22)(43,56)
\qbezier(131,56)(43,56)(43,56)
\qbezier(131,56)(95,22)(95,22)

\put(66,20){\vector(0,-1){9}}
\multiput(66,38)(0,-4){5}{\line(0,-2){2}}
\put(74,10){$x$}

\end{picture}
\end{minipage}
\hspace{1cm}$=$\hspace{1cm}
\begin{minipage}{130pt}
\begin{picture}(130,70)

\qbezier(7,22)(7,22)(95,22)
\qbezier(7,22)(7,22)(43,56)
\qbezier(131,56)(43,56)(43,56)
\qbezier(131,56)(95,22)(95,22)

\put(69,38){\vector(0,1){26}}
\put(77,65){$\rho(x)$}

\end{picture}
\end{minipage}

\end{center}
\caption{A basic inversion}
\label{Fig.B}
\end{figure}

A {\em basic inversion} is an operation which reverses the normal orientation of a semi-sheet and changes the element $x$ assigned the sheet by $\rho(x)$.  
See Figure~\ref{Fig.B}.   The coloring conditions are preserved under basis inversions.  

An {\em $(X, \rho)$-coloring} of $D$ is the equivalence class of 
an assignment of a normal orientation and an element of $X$ to each semi-sheet of $D$ satisfying the coloring conditions.  
Here the equivalence relation is generated by basic inversions.    

\begin{figure}[H]
\begin{center}
\begin{minipage}{160pt}
\begin{picture}(160,80)

\qbezier(7,32)(7,32)(95,32)
\qbezier(7,32)(7,32)(43,66)
\qbezier(131,66)(43,66)(43,66)
\qbezier(131,66)(95,32)(95,32)

\put(66,30){\vector(0,-1){9}}
\multiput(66,48)(0,-4){5}{\line(0,-2){2}}
\put(74,20){$x$}
\put(60,75){$[y_1]$}
\put(60,5){$[y_2]$}
\put(140,35){\framebox(46,14){$y_1^{x} = y_2$} }

\end{picture}
\end{minipage}
\end{center}
\caption{Coloring condition for regions}
\label{Fig.D}
\end{figure}
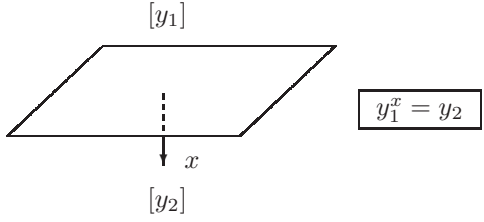

Let $Y$ be an $(X, \rho)$-set. An {\em $(X, \rho)_Y$-coloring} of $D$ is an 
$(X, \rho)$-coloring with an assignment of an element of $Y$ to each complementary region of $D$ 
satisfying the following condition (Figure~\ref{Fig.D}): 
\begin{itemize}
\item Suppose that two adjacent regions $f_1$ and $f_2$ which are separated by a semi-sheet, say $e$, are labeled by $y_1$ and $y_2$.  Suppose that the semi-sheet $e$ is labeled by $x$.  If the normal orientation of $e$ points from $f_1$ to $f_2$, then $y_1^x = y_2$. 
\end{itemize}

In figures, the label $y \in Y$ assigned a region is indicated by $[y]$. See Figure~\ref{Fig.D}. 

\begin{proposition} 
Let $(X,\rho)$ be a symmetric quandle, and $Y$  an $(X,\rho)$-set. 
If two diagrams represent the same unoriented surface-link type, then there is a bijection between the sets of $(X,\rho)$-colorings of the diagrams, and there is a bijection between the sets of $(X, \rho)_Y$-colorings.
\end{proposition}

\begin{proof}
Suppose that $D$ and $D'$ are diagrams related by a single Roseman move \cite{CarterSaito, Roseman}. 
Let $E$ be a $3$-disk in $\mathbb{R}^3$ in which the Roseman move is applied.
For each $(X,\rho)_Y$-coloring of $D$, its restriction to $D\setminus E (=D'\setminus E)$ 
can be uniquely extended to an $(X,\rho)_Y$-coloring of $D'$.
Thus there is a bijection between the set of $(X,\rho)_Y$-colorings of $D$ and that of $D'$. 
\end{proof}

Let $D$ be an unoriented surface-link diagram.  
Fix an $(X, \rho)_Y$-coloring of $D$, say $C$. 
For a triple point $v$ of $D$, there are eight complementary regions of $D$ around $v$. 
(Some of them may be the same.)
Choose one of them, say $f$, which we call {\it a specified region} for $v$,  and let $y$ be the label of $f$.  

Let $e_1$, $e_2$, and $e_3$ be the lower semi-sheet, the middle semi-sheet, and the upper semi-sheet at $v$, respectively,  which face the region $f$.  
By basic inversions, we assume that the normal orientations $n_1$, $n_2$ and $n_3$  of them point from $f$ to the opposite regions.  
Let $x_1$, $x_2$  and $x_3$ be the labels of them, respectively.  
The {\em sign} of $v$ 
 with respect to the region $f$ is $+1$ (or $-1$)  if the pair of normal orientations $(n_3, n_2, n_1)$ does (or does not) match the orientation of ${\mathbb R}^3$.  
The {\em weight} of $v$ is defined to be $\epsilon (y, x_1, x_2, x_3)$ where $\epsilon $ is the sign of $v$. See Figure~\ref{Fig.E}.

\begin{figure}[H]
\begin{center}
\begin{minipage}{140pt}
\begin{picture}(130,140)

\qbezier(7,57)(7,57)(95,57)
\qbezier(7,57)(7,57)(29,78)
\qbezier(131,91)(95,57)(95,57)
\qbezier(131,91)(131,91)(110,91)

\qbezier(50,58)(50,58)(50,100)
\qbezier(84,132)(50,100)(50,100)
\qbezier(84,132)(84,132)(84,117)
\qbezier(50,58)(50,58)(69,76)

\qbezier(50,18)(50,18)(50,53)
\qbezier(50,18)(50,18)(69,36)
\qbezier(54,57)(50,53)(50,53)

\qbezier(71,78)(71,78)(110,78)
\qbezier(71,78)(71,78)(71,117)
\qbezier(71,117)(71,117)(110,117)
\qbezier(110,78)(110,78)(110,117)

\qbezier(71,38)(71,38)(110,38)
\qbezier(71,38)(71,38)(71,57)
\qbezier(110,38)(110,38)(110,71)

\qbezier(30,38)(30,38)(50,38)
\qbezier(30,38)(30,38)(30,57)

\qbezier(30,78)(30,78)(30,117)
\qbezier(30,117)(30,117)(66,117)
\qbezier(30,78)(30,78)(50,78)
\qbezier(66,116)(66,117)(66,117)

\put(35,-2){\framebox(60,14){$(y,x_1,x_2,x_3)$}}

\put(19,44){\colorbox{white}{\textcolor{black}{$x_1$}} }
\put(28,28){$[y]$}
\put(60,39){$x_2$}
\put(27,70){$x_3$}

\put(55,47){\vector(1,0){25}}
\put(37,65){\vector(0,1){25}}
\qbezier(35,45)(35,45)(35.5,45.5)
\qbezier(37,47)(37,47)(37.5,47.5)
\qbezier(39,49)(39,49)(39.5,49.5)
\put(44,54){\vector(1,1){0}}

\end{picture}
\end{minipage}
\hspace{1cm}
\begin{minipage}{140pt}
\begin{picture}(130,140)

\qbezier(7,57)(7,57)(95,57)
\qbezier(7,57)(7,57)(29,78)
\qbezier(131,91)(95,57)(95,57)
\qbezier(131,91)(131,91)(110,91)

\qbezier(50,58)(50,58)(50,100)
\qbezier(84,132)(50,100)(50,100)
\qbezier(84,132)(84,132)(84,117)
\qbezier(50,58)(50,58)(69,76)

\qbezier(50,18)(50,18)(50,53)
\qbezier(50,18)(50,18)(69,36)
\qbezier(54,57)(50,53)(50,53)

\qbezier(71,78)(71,78)(110,78)
\qbezier(71,78)(71,78)(71,117)
\qbezier(71,117)(71,117)(110,117)
\qbezier(110,78)(110,78)(110,117)

\qbezier(71,38)(71,38)(110,38)
\qbezier(71,38)(71,38)(71,57)
\qbezier(110,38)(110,38)(110,71)

\qbezier(30,38)(30,38)(50,38)
\qbezier(30,38)(30,38)(30,57)

\qbezier(30,78)(30,78)(30,117)
\qbezier(30,117)(30,117)(66,117)
\qbezier(30,78)(30,78)(50,78)
\qbezier(66,116)(66,117)(66,117)

\put(55,90){\vector(1,0){25}}
\put(54,117){\vector(1,1){10}}

\qbezier(42,105)(42,105)(44,107)
\qbezier(46,109)(46,109)(48,111)
\qbezier(50,113)(50,113)(52,115)
\put(35,57){\vector(0,-1){10}}
\multiput(35,60)(0,4){3}{\line(0,-1){2.0}}
\put(50,130){$x_1$}
\put(80,95){$x_2$}
\put(18,44){$x_3$}
\put(34,87){$[y]$}
\put(36,-2){\framebox(65,14){$-(y,x_1,x_2,x_3)$}}

\end{picture}
\end{minipage}
\end{center}
\caption{Weights}
\label{Fig.E}
\end{figure}

\begin{lemma}
As an element of $C_3^{{\rm R}, \rho}(X)_Y = C_3(X)_Y / D_3^{\rho}(X)_Y$, the weight of $v$ does not depend on the specified region. 
\end{lemma}

\begin{proof}
When we change the specified region, the difference of the weights belongs to $D_3^{\rho}(X)_Y$.\end{proof}

For a diagram $D$ and an $(X,\rho)_Y$-coloring $C$, we define a chain $c_{D, C}$ by 
$$
 c_{D, C} = \sum_{v} \epsilon (y, x_1, x_2, x_3) \quad \in C_3^{{\rm R}, \rho}(X)_Y \  {\rm or}\ C_3^{{\rm Q},\rho}(X)_Y, 
$$
where $v$ runs over all triple points of $v$ and $\epsilon (y, x_1, x_2, x_3)$ is the weight of $v$. 

Let $(D,C)$ be a pair of an unoriented surface-link diagram and an $(X,\rho)_Y$-coloring.  If $D$ changes into $D'$ by a single Roseman move in a $3$-disk support $E$, then there is a unique $(X,\rho)_Y$-coloring, say $C'$, of $D'$ which is identical with $C$ outside $E$.   In this situation, we say that $(D',C')$ is obtained from $(D,C)$ by a Roseman move (with support $E$).

We say that $(D, C)$ and $(D', C')$ are {\em Roseman move equivalent} if they are related by a finite sequence of Roseman moves.  The equivalence class of $(D, C)$ is called an {\em $(X, \rho)_Y$-colored surface-link}.  

\begin{theorem}
\label{theorem:3-cycle} 
The chain $c_{D,C}$ is a $3$-cycle of $C_\ast^{{\rm Q},\rho}(X)_Y$, and 
the homology class $[c_{D,C}] \in H_3^{{\rm Q}, \rho}(X)_Y$ is an invariant of an $(X, \rho)_Y$-colored surface-link.  
\end{theorem}
 
\begin{proof} 
First we show that the chain $c_{D,C}$ is a $3$-cycle of  $C_\ast^{{\rm Q}, \rho}(X)_Y$.  
Let us fix a checkerboard coloring of the complementary regions of $D$. 
Using basic inversions, we assume that normal orientations of semi-sheets of $D$ point from the black regions to the white.
Give an orientation to each double curve of $D$ such that 
the orientation vector followed by the normal orientation vector of the over semi-sheet and the normal orientation vector of the under semi-sheet  matches the orientation of ${\mathbb R}^3$.  See Figure~\ref{Fig.G}.

For a double curve $\gamma$, choose a black region around $\gamma$, and let $y$ be the coloring assigned to this region. 
Let $x_1$ and $x_2$ be the color of the under semi-sheet and the color of the over semi-sheet around $\gamma$ which face the black region.  
Assign a triple $(y, x_1, x_2)\in Y\times X^2$ to the double curve $\gamma$. See Figure~\ref{Fig.G}.  It depends on a choice of the black region around $\gamma$, however it determines a unique element of  $C_\ast^{{\rm Q}, \rho}(X)_Y$.

\begin{figure}[H]
\begin{center}
\begin{minipage}{160pt}
\begin{picture}(160,130)

\put(12,70){\setlength{\fboxrule}{17pt}\fcolorbox[gray]{.9}{.9}{\hspace{-5pt}}}
\put(61,30){\setlength{\fboxrule}{12pt}\fcolorbox[gray]{.9}{.9}{\ \ \ }}

\qbezier(7,47)(7,47)(43,81)
\qbezier(7,47)(7,47)(46,47)
\qbezier(46,47)(46,47)(50,51)
\qbezier(50,81)(43,81)(43,81)

\qbezier(131,81)(131,81)(95,47)
\qbezier(53,47)(53,47)(95,47)
\qbezier(53,47)(53,47)(89,81)
\qbezier(131,81)(89,81)(89,81)

\qbezier(50,8)(50,8)(50,90)
\qbezier(84,122)(50,90)(50,90)
\qbezier(50,8)(50,8)(84,40)
\qbezier(84,47)(84,40)(84,40)
\qbezier(84,122)(84,78)(84,78)

\put(65,88){\vector(1,0){25}}
\put(35,52){\vector(0,-1){15}}
\multiput(35,63)(0,-4){3}{\line(0,-2){2}}
\put(25,30){$x_1$}
\put(94,89){$x_2$}
\put(23,75){$[y]$}

\thicklines
\qbezier(51,50)(51,50)(57,56)
\qbezier(60,58.5)(60,58.5)(66,64.5)
\qbezier(69,67)(69,67)(75,72.5)
\qbezier(78,75.5)(78,75.5)(83.5,80)
\put(84,80.5){\vector(1,1){0}}
\thinlines

\put(140,65){\ovalbox{$(y,x_1,x_2)$}}

\qbezier(135,67)(100,72)(79,70)
\put(79,70){\vector(-1,0){0}}

\end{picture}
\end{minipage}
\end{center}
\caption{The label of a double point set}
\label{Fig.G}
\end{figure}
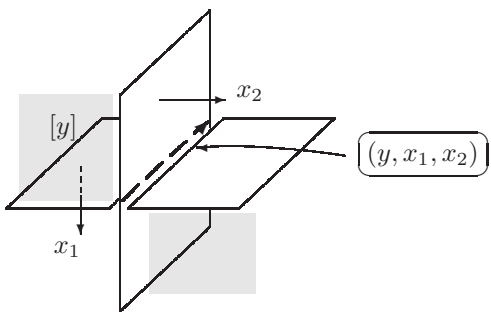

If the double curve is not a simple loop, then we give the initial point the weight $-(y, x_1, x_2)$ and the terminal point the weight $(y, x_1, x_2)$.
For each triple point $v$ of $D$, we choose its specified region from the black regions.
In the situation of Figure~\ref{Fig.I},  the weights of $v$ is $(y,x_1,x_2, x_3)$, and
$$
\begin{array}{lll}
\partial _3(y,x_1,x_2, x_3)
&=&
 -(y,x_2, x_3)
 +(y^{x_1},x_2, x_3) 
 +(y,x_1, x_3 )
 -(y^{x_2},x_1^{x_2}, x_3 ) \\
&  &-(y, x_1, x_2) 
 +(y^{x_3}, x_1^{x_3}, x_2^{x_3}) \\
&\equiv & 
 -(y,x_2, x_3)
 -(y^{x_1 x_3},x_2^{x_3}, \rho(x_3)) 
 +(y,x_1, x_3 )
 +(y^{x_1 x_2},\rho(x_1)^{x_2}, x_3 ) \\
& & -(y, x_1, x_2) 
 -(y^{x_1 x_3}, \rho(x_1)^{x_3}, x_2^{x_3}), 
 \end{array}
$$
mod $D_2^{\rho}(X)_Y.$
The six terms are exactly the same with the weights assigned the endpoints of double curves around $v$.  
Since the weights of endpoints of double curves connecting triple points 
must be canceled, we see that $\partial _3(c_{D,C}) \in C_2^{{\rm Q},\rho}(X)_Y$ is the sum of the weights of initial or terminal points of 
the double curves one of whose endpoints are triples points and the other endpoints are branch points. However such weights are in the form $\pm (y, x, x)$  or $\pm (y, x, \rho(x))$ which vanish in $C_2^{{\rm Q},\rho}(X)_Y$.  Therefore the chain $c_{D,C}$ is a $3$-cycle of $C_\ast^{{\rm Q},\rho}(X)_Y$.

\begin{figure}[H]
\begin{center}
\begin{minipage}{130pt}
\begin{picture}(130,130)

\put(60,75){\setlength{\fboxrule}{18pt}\fcolorbox[gray]{.9}{.9}{}}
\put(18,28){\setlength{\fboxrule}{13pt}\fcolorbox[gray]{.9}{.9}{}}
\thicklines
\multiput(66,108)(-1,0){36}{\color[gray]{0.9}{\line(1,1){11}}}
\multiput(113,30)(0.5,0.5){37}{\color[gray]{0.9}{\line(0,1){30}}}
\thinlines

\qbezier(7,47)(7,47)(95,47)
\qbezier(7,47)(7,47)(29,68)
\qbezier(131,81)(95,47)(95,47)
\qbezier(131,81)(131,81)(110,81)

\qbezier(50,48)(50,48)(50,90)
\qbezier(84,122)(50,90)(50,90)
\qbezier(84,122)(84,122)(84,107)
\qbezier(50,48)(50,48)(69,66)

\qbezier(50,8)(50,8)(50,43)
\qbezier(50,8)(50,8)(69,26)
\qbezier(54,47)(50,43)(50,43)

\qbezier(71,68)(71,68)(110,68)
\qbezier(71,68)(71,68)(71,107)
\qbezier(71,107)(71,107)(110,107)
\qbezier(110,68)(110,68)(110,107)

\qbezier(71,28)(71,28)(110,28)
\qbezier(71,28)(71,28)(71,47)
\qbezier(110,28)(110,28)(110,61)

\qbezier(30,28)(30,28)(50,28)
\qbezier(30,28)(30,28)(30,47)

\qbezier(30,68)(30,68)(30,107)
\qbezier(30,107)(30,107)(66,107)
\qbezier(30,68)(30,68)(50,68)
\qbezier(66,106)(66,107)(66,107)

\put(55,37){\vector(1,0){20}}
\put(37,55){\vector(0,1){20}}
\qbezier(35,35)(35,35)(35.5,35.5)
\qbezier(37,37)(37,37)(37.5,37.5)
\qbezier(39,39)(39,39)(39.5,39.5)
\put(44,44){\vector(1,1){0}}

\put(15,84){\footnotesize\colorbox {white}{\textcolor{black}{$\rho(x_2)^{x_3}$}} }

\put(50,80){\vector(-1,0){8}}
\multiput(55,80)(4,0){3}{\line(-1,0){2.0}}

\put(99,107){\vector(1,1){5}}

\qbezier(87,95)(87,95)(89,97)
\qbezier(91,99)(91,99)(93,101)
\qbezier(95,103)(95,103)(97,105)
\qbezier(99,107)(99,107)(101,109)

\put(80,47){\vector(0,-1){8}}
\multiput(80,50)(0,4){3}{\line(0,-1){2.0}}

\put(19,34){\colorbox [gray]{0.9}{\textcolor{black}{$x_1$}} }
\put(28,18){$[y]$}
\put(60,29){$x_2$}
\put(27,60){$x_3$}
\put(85,39){\small$\rho(x_3)$}
\put(100,118){\footnotesize$x_1^{x_2 x_3}$}
\put(75,80){\small$[y^{x_2 x_3}]$}
\put(39,112){\small$[y^{x_1 x_3}]$}
\put(111,50){\tiny$[y^{x_1 x_2}]$}

\end{picture}
\end{minipage}
\hspace{1cm}
\begin{minipage}{200pt}

\begin{picture}(200,130)

\put(90,20){\line(0,1){90}}
\put(45,65){\line(1,0){90}}
\put(90,65){\line(1,1){33}}
\put(90,65){\line(-1,-1){33}}

\put(90,50){\vector(0,1){0}}
\put(90,80){\vector(0,-1){0}}
\put(80,55){\vector(1,1){0}}
\put(100,75){\vector(-1,-1){0}}
\put(60,65){\vector(-1,0){0}}
\put(120,65){\vector(1,0){0}}

\put(70,2){\ovalbox{$-(y,x_1,x_2)$}}
\put(4,14){\ovalbox{$-(y,x_2,x_3)$}}
\put(-10,60){\ovalbox{$(y,x_1,x_3)$}}
\put(40,120){\ovalbox{$-(y^{x_1 x_3},\rho(x_1)^{x_3},x_2^{x_3})$}}
\put(130,100){\ovalbox{$-(y^{x_1 x_3},x_2^{x_3},\rho(x_3))$}}
\put(140,60){\ovalbox{$(y^{x_1 x_2}, \rho(x_1)^{x_2}, x_3)$}}
\end{picture}
\end{minipage}
\end{center}
\caption{The label of a double curve}
\label{Fig.I}
\end{figure}
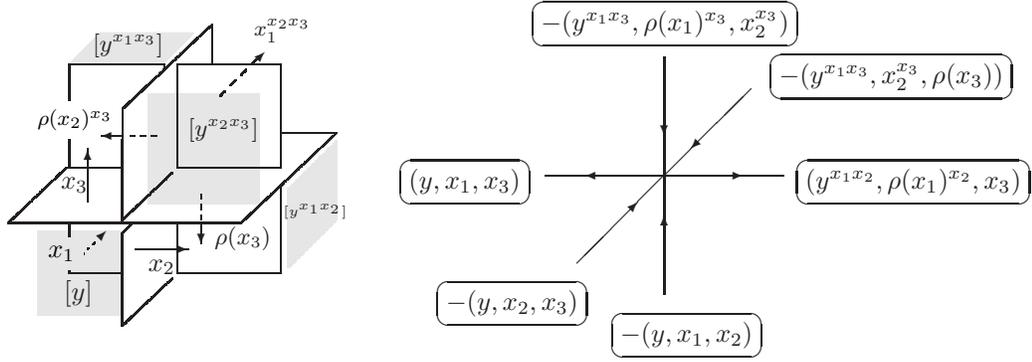

We show that if $(D, C)$ changes into $(D', C')$ by a Roseman move (whose support is a $3$-disk $E$),  then 
$[c_{D, C}]= [c_{D', C'}]$ in  $H_3^{{\rm Q}, \rho}(X)_Y$. 
The diagrams $D$ and $D'$ and their colorings $C$ and $C'$ are identical outside of $E$.   The preimage of $E$ in the surface-link $F$ and $F'$ described by $D$ and $D'$ are disks in $F$ and $F'$, respectively.  Assign the disks orientations and project the orientation to the semi-sheets of $D$ and $D'$  restricted in $E$.   Assume that  the normal orientation of semi-sheets which meets $E$ are compatible with these orientations with respect to the orientation of $\mathbb{R}^3$.  For each triple point $v$ in $E$, there is a unique complementary region around $v$  such that the semi-sheets face the region have normal orientations from the region.  When we choose such a region as the specified region, the sum of weights of triple points of $D$ in $E$ and that of triple points of $D'$ in $E$ differ by an element of $D_3^{\rm Q}(X)$ or $\partial(C_4^{\rm Q}(X))$ by the same reason as in \cite{CJKLS03, CKS01}.  Thus the difference of the weight sum for $(D, C)$ and for $(D', C')$ belongs to 
$D_3^{\rm Q}(X) + D_3^{\rho }(X) + \partial(C_4^{\rm Q}(X))$.  
\end{proof}

Let 
$$
{\cal H}(D) = \{ [c_{D, C}] \in H_3^{{\rm Q}, \rho}(X)_Y \, 
| \, C: \mbox{$(X, \rho)_Y$-colorings of $D$} \} 
$$
as a multi-set.    
For  a $3$-cocycle $\theta$ of the cochain complex $C^*_{{\rm Q}, \rho}(X, A)_Y$ with a coefficient group $A$, let 
$$
\Phi _\theta(D) = \{ \theta(c_{D, C}) \in A \, 
| \, C: \mbox{$(X, \rho)_Y$-colorings of $D$} \} 
$$
as a multi-set.  

\begin{corollary}\label{theorem:invariant2dim}
The multi-sets ${\cal H}(D)$ and $\Phi _\theta(D)$ are invariants of the surface-link type of $D$.  
\end{corollary}

We recall quandle cocycle invariants of oriented surface-links in the sense of 
\cite{CJKLS03, CKS01, RS1}. 
Let $D^+$ be an oriented diagram such that $D$ is obtained from $D^+$ by forgetting orientation.  Using the orientation, we assign 
the semi-sheets of $D^+$ normal orientations such that the orientation vectors followed by the normal orientation vector matches the orientation of ${\mathbb R}^3$.  An $X_Y$-coloring is assignment of an element of $X$ to each semi-sheet and an element of $Y$ to each complementary region of $D$ satisfying the conditions illustrated in 
the left of Figure~\ref{Fig.C} and in Figure~\ref{Fig.D}.   The weight of each triple point is defined as in Figure~\ref{Fig.E}.  Let 
$$
{\cal H}^{\rm ori}(D^+) = \{ [c_{D^+, C}] \in H_3^{{\rm Q}}(X)_Y \, 
| \, C: \mbox{$X_Y$-colorings of $D^+$} \} 
$$
and 
$$
\Phi _\theta^{\rm ori}(D^+) = \{ \theta(c_{D^+, C}) \in A \, 
| \, C: \mbox{$X_Y$-colorings of $D^+$} \} 
$$
as multi-sets, where $\theta$ is a $3$-cocycle of $C^*_{{\rm Q}}(X, A)_Y$.  

If $\theta:\mathbb{Z}(Y\times X^3)\to A$ is a $3$-cocycle of $C_{{\rm Q},\rho}^*(X, A)_Y$, then $\theta$ is a $3$-cocycle of $C_{{\rm Q}}^*(X, A)_Y$. 

\begin{theorem}\label{thm:equality2dim}
Let $(X, \rho)$ be a symmetric quandle, $Y$ an $(X, \rho)$-set, $A$ an abelian group, and $\theta$ a $3$-cocycle of $C_{{\rm Q},\rho}^*(X;A)_Y$.  Let $D$ be a diagram of an unoriented surface-link $F$  such that every component of $F$ is orientable.  Let $D^+$ be an arbitrarily oriented diagram of $D$.  
Then $$ 
\Phi _\theta(D)= \Phi^{\rm ori}_\theta(D^+).$$
\end{theorem}

\begin{proof}
Assign a normal orientation to each semi-sheet of $D$ (and $D^+$) such that the orientation vectors determined from $D^+$ followed by the normal orientation vector match the orientation of ${\mathbb R}^3$.  For 
each  $(X, \rho)_Y$-coloring $C$ of $D$, take a representative of $C$ such that the normal orientation of each semi-sheet of $D$ is the same with the normal orientation.   Then $C$ determines an assignment of an element of $X$ to each semi-sheet.   Let $C^+$ be this assignment together with the same assignment of an element of $Y$ to each   complementary region as $C$.  Then $C^+$  is an $X_Y$-coloring for the oriented diagram $D^+$ in the sense of \cite{CJKLS03, CKS01, RS1}.  Then there is a bijection between the set of $(X, \rho)_Y$-colorings of $D$ and the set of  
$X_Y$-coloring of $D^+$.  At each triple point the weights for $(D,C)$ and $(D', C')$ evaluated by $\theta$ are the same,  and we have 
$$ \theta(c_{D, C}) = \theta(c_{D^+, C^+}). $$
and  $ 
\Phi_\theta(D^+) =  \Phi^{\rm ori} _\theta(D)$. 
\end{proof}

\begin{theorem}\label{theorem:triplepoint}
Let $(X, \rho)$ be a symmetric quandle, $Y$  an $(X, \rho)$-set and $\theta$ a  $3$-cocycle of $C_{{\rm Q}, \rho}^\ast(X, {\mathbb Z})_Y$.  
Suppose that for any $(y, x_1, x_2, x_3) \in Y \times X^3$, $|\theta( (e, x_{i_1}, x_{i_2}, x_{i_3}) )| \leq 1$.  
Let  $F$ be a surface-link and $D$ a diagram of $F$ 
with an $(X, \rho)_Y$-coloring $C$.  
Then  $| \theta (c_{D,C}) | \leq t(F)$, 
where $t(F)$ is the minimal triple point number of $F$.   
\end{theorem} 

\begin{proof} 
Let $D'$ be a diagram of $F$ such that the number of triple points is $t(F)$.  Let $C'$ be an $(X, \rho)_Y$-coloring of $D'$ such that $(D', C')$ is Roseman move equivalent to $(D,C)$.   By Theorem~\ref{theorem:3-cycle}, $\theta(c_{D,C}) = \theta(c_{D', C'})$.  
By definition, the chain  $c(D', C')$ is expressed as $\sum_{i=1}^{t(F)} \epsilon_i (y_i, x_{i_1} x_{i_2}, x_{i_3})$ for some 
elements $(y_i, x_{i_1} x_{i_2}, x_{i_3})$ of $Y \times X^3$ and $\epsilon_i \in \{\pm1\}$.  
Since $|\theta( (y_i, x_{i_1} x_{i_2}, x_{i_3}) )| \leq 1$, we have $| \theta(c_{D',C'})| \leq t(F)$.  
\end{proof} 

This theorem can be used to estimate the minimal triple point numbers of surface-links.  For example, see Proposition~\ref{proposition:cocycledihedral4},  Theorems~\ref{theorem:triple} and \ref{theorem:P2P2}. 

\section{Examples and applications}
\label{Examplescocycles}

First,  we give some examples of quandle cocycles. 

\begin{example}[\cite{Oshiro}]\label{example:linking2dim}{\rm 
Let $X$ be the order $6$ trivial quandle $T_6=  \{e_1,e'_1,e_2,e'_2, e_3, e'_3\}$,  and  
$\rho: X \to X$ a good involution with $\rho(e_i) = e'_i$ $(i \in \{1,2,3\})$.  Let $Y=\{e\}$, which is an $(X, \rho)$-set. 
Define a map $\theta:Y\times X^3 \to {\mathbb Z}$ by 
$$  
\begin{array}{ll} 
\theta 
= & \chi_{(e, e_1, e_2, e_3)} + \chi_{(e, e'_1, e'_2, e_3)}   + \chi_{(e, e'_1, e_2, e'_3)} + \chi_{(e, e_1, e'_2, e'_3)} \\
& - \chi_{(e, e'_1, e_2, e_3)} - \chi_{(e, e_1, e'_2, e_3)} - \chi_{(e, e_1, e_2, e'_3)} - \chi_{(e, e'_1, e'_2, e'_3)}
\end{array}
$$
where $\chi_{(e, a, b,c)}$ is defined by 
$\chi_{(e,a,b,c)} (e,x,y,z) = 1$ if $(e, x ,y,z)=(e, a, b,c)$ and $\chi_{(e,a,b,c)} (e,x,y,z) = 0$ otherwise. 
The linear extension $\theta: \mathbb{Z}(Y\times X^3)\to {\mathbb Z}$  is a $3$-cocycle.   In \cite{Oshiro}  it is shown that the invariant $\Phi_\theta$ has information of the triple linking invariants \cite{CKSS01} for oriented surface-links.   
}\end{example}

\begin{example}[\cite{Oshiro}]\label{example:mod2linking2dim}{\rm 
Let $X$ be the order $2$ trivial quandle $T_2=  \{e_1,e_2\}$,  and  
$\rho: X \to X$ a good involution with $\rho(x) = x$ $(x \in X)$.  Let $Y=\{e\}$, which is an $(X, \rho)$-set. 
Define a map $\theta:Y\times X^3 \to {\mathbb Z}/ 2{\mathbb Z}$ by 
$$  
\theta 
=  \chi_{(e, e_1, e_2, e_1)}.     
$$
The linear extension  $\theta: \mathbb{Z}(Y\times X^3)\to {\mathbb Z}/ 2{\mathbb Z}$ is a $3$-cocycle.   In \cite{Oshiro} it is shown that the invariant $\Phi_\theta$ has information of the mod $2$ triple linking invariants \cite{Satoh2002} for surface-links which are not necessarily oriented or orientable.  
}\end{example}

\begin{example}\label{example:cocycledihedral4}{\rm 
Let $X$ be the order $4$ dihedral quandle $R_4$, which we rename the elements $0,1,2,3$ by $e_1, e_2, e'_1, e'_2$ respectively.  
Let $\rho: X \to X$ be the antipodal map, i.e., $\rho(e_i)=e'_i$ $(i=1,2)$.   Let $Y=\{e\}$, which is an $(X, \rho)$-set. 
Define a map $\theta:Y\times X^3 \to {\mathbb Z}$ by 
$$  
\begin{array}{ll} 
\theta 
= 
& +\chi_{(e, e_1, e_2, e_1)} + \chi_{(e, e'_1, e'_2, e_1)}   + \chi_{(e, e'_1, e_2, e'_1)} + \chi_{(e, e_1, e'_2, e'_1)} \\ 
& - \chi_{(e, e'_1, e_2, e_1)} - \chi_{(e, e_1, e'_2, e_1)} - \chi_{(e, e_1, e_2, e'_1)}  - \chi_{(e, e'_1, e'_2, e'_1)} \\ 
& - \chi_{(e, e_2, e_1, e_2)} - \chi_{(e, e'_2, e'_1, e_2)}   - \chi_{(e, e'_2, e_1, e'_2)} - \chi_{(e, e_2, e'_1, e'_2)} \\
& + \chi_{(e, e'_2, e_1, e_2)} + \chi_{(e, e_2, e'_1, e_2)} + \chi_{(e, e_2, e_1, e'_2)}  + \chi_{(e, e'_2, e'_1, e'_2)}. 
\end{array}
$$
The linear extension  $\theta: \mathbb{Z}(Y\times X^3)\to {\mathbb Z}$ is a $3$-cocycle.  
}\end{example}

Now we give some applications on minimal triple point numbers of surface-links. The following is obtained from Theorem~\ref{theorem:triplepoint}.  

\begin{proposition}\label{proposition:cocycledihedral4}
Let $(X, \rho)$, $Y$, and $\theta$ be the symmetric quandle, the $(X, \rho)$-set, and the $3$-cocycle given in Example~$\ref{example:cocycledihedral4}$.   
Suppose that a surface-link $F$ has a diagram $D$ with an $(X, \rho)_Y$-coloring $C$.  
Then $| \theta (c_{D,C}) | \leq t(F)$, 
where $t(F)$ is the minimal triple point number of $F$.   
\end{proposition}

\begin{figure}[H]
\begin{center}

\begin{minipage}{50pt}
\begin{picture}(40,60)
\put(6,6){\line(0,1){48}}

\qbezier(6,6)(18,10)(18,30)
\qbezier(6,54)(18,50)(18,30)

\put(34,6){\line(0,1){48}}

\qbezier(34,6)(22,10)(22,30)
\qbezier(34,54)(22,50)(22,30)

\put(3,30){\vector(1,0){10}}
\put(37,30){\vector(-1,0){10}}
\put(0,33){\small $x$}
\put(36,33){\small $y$}

\end{picture}
\end{minipage}
\hspace{-10pt}
\begin{minipage}{60pt}
\begin{picture}(50,60)
\put(6,6){\line(0,1){48}}
\put(6,6){\line(1,0){14}}
\put(6,54){\line(1,0){14}}
\put(18,22){\line(0,1){16}}

\qbezier(20,54)(23.5,53.7)(25,52.5)

\put(18,38){\line(1,0){2}}
\qbezier(20,38)(23.2,38.5)(24,39)

\qbezier(25,52.5)(28,51)(28,46)
\qbezier(28,46)(27.8,41.5)(26,40.5)

\put(18,22){\line(1,0){2}}
\qbezier(20,22)(23.5,21.7)(24,21)

\qbezier(20,6)(23.5,6.5)(25,7.5)

\qbezier(26.7,19.3)(28,18)(28,14)
\qbezier(28,14)(28,9)(25,7.5)

\put(44,6){\line(0,1){48}}
\put(44,6){\line(-1,0){14}}
\put(44,54){\line(-1,0){14}}
\put(32,22){\line(0,1){16}}

\qbezier(30,54)(26.5,53.7)(26.5,53.3)

\put(32,38){\line(-1,0){2}}
\qbezier(30,38)(26.5,38.5)(25,39.5)

\qbezier(23.8,52)(22,50)(22,46)
\qbezier(22,46)(22,41)(25,39.5)

\qbezier(30,22)(26.5,21.7)(25,20.5)

\put(32,22){\line(-1,0){2}}
\qbezier(30,6)(26.5,6.5)(26.5,7)

\qbezier(25,20.5)(22,19)(22,14)
\qbezier(22,14)(22.7,9.3)(23.6,9)

\end{picture}
\end{minipage}
\hspace{-7pt}
\begin{minipage}{60pt}
\begin{picture}(50,60)
\put(6,6){\line(0,1){21}}
\put(6,33){\line(0,1){21}}
\qbezier(6,33)(6,32)(8,30)
\qbezier(8,30)(13,23.5)(14,30)
\qbezier(14,30)(12,34)(9,31)
\qbezier(6,27)(6,28)(7,29)

\put(6,6){\line(1,0){14}}
\put(6,54){\line(1,0){14}}
\put(18,22){\line(0,1){16}}

\qbezier(20,54)(23.5,53.7)(25,52.5)

\put(18,38){\line(1,0){2}}
\qbezier(20,38)(23.2,38.5)(24,39)

\qbezier(25,52.5)(28,51)(28,46)
\qbezier(28,46)(27.8,41.5)(26,40.5)

\put(18,22){\line(1,0){2}}
\qbezier(20,22)(23.5,21.7)(24,21)

\qbezier(20,6)(23.5,6.5)(25,7.5)

\qbezier(26.7,19.3)(28,18)(28,14)
\qbezier(28,14)(28,9)(25,7.5)

\put(44,6){\line(0,1){21}}
\put(44,33){\line(0,1){21}}
\qbezier(43.8,33)(44,32)(42,30)
\qbezier(42,30)(37,23.5)(36,30)
\qbezier(36,30)(38,34)(41,31)
\qbezier(43.8,27)(44,28)(43,29)

\put(44,6){\line(-1,0){14}}
\put(44,54){\line(-1,0){14}}
\put(32,22){\line(0,1){16}}

\qbezier(30,54)(26.5,53.7)(26.5,53.3)

\put(32,38){\line(-1,0){2}}
\qbezier(30,38)(26.5,38.5)(25,39.5)

\qbezier(23.8,52)(22,50)(22,46)
\qbezier(22,46)(22,41)(25,39.5)

\qbezier(30,22)(26.5,21.7)(25,20.5)

\put(32,22){\line(-1,0){2}}
\qbezier(30,6)(26.5,6.5)(26.5,7)

\qbezier(25,20.5)(22,19)(22,14)
\qbezier(22,14)(22.7,9.3)(23.6,9)

\end{picture}
\end{minipage}
\hspace{-7pt}
\begin{minipage}{60pt}
\begin{picture}(50,60)
\put(6,6){\line(0,1){20}}
\put(6,34){\line(0,1){20}}

\put(6,6){\line(1,0){14}}
\put(6,54){\line(1,0){14}}

\put(18,22){\line(0,1){4}}
\put(18,34){\line(0,1){4}}

\qbezier(6,34)(6,34)(18,26)
\qbezier(6,26)(6,26)(10,28.7)
\qbezier(14,31)(18,34)(18,34)

\qbezier(20,54)(23.5,53.7)(25,52.5)

\put(18,38){\line(1,0){2}}
\qbezier(20,38)(23.2,38.5)(24,39)

\qbezier(25,52.5)(28,51)(28,46)
\qbezier(28,46)(27.8,41.5)(26,40.5)

\put(18,22){\line(1,0){2}}
\qbezier(20,22)(23.5,21.7)(24,21)

\qbezier(20,6)(23.5,6.5)(25,7.5)

\qbezier(26.7,19.3)(28,18)(28,14)
\qbezier(28,14)(28,9)(25,7.5)

\put(44,6){\line(0,1){20}}
\put(44,34){\line(0,1){20}}

\put(44,6){\line(-1,0){14}}
\put(44,54){\line(-1,0){14}}

\put(32,22){\line(0,1){4}}
\put(32,34){\line(0,1){4}}

\qbezier(44,34)(44,34)(32,26)
\qbezier(44,26)(44,26)(40,28.7)
\qbezier(36,31)(32,34)(32,34)

\qbezier(30,54)(26.5,53.7)(26.5,53.3)

\put(32,38){\line(-1,0){2}}
\qbezier(30,38)(26.5,38.5)(25,39.5)

\qbezier(23.8,52)(22,50)(22,46)
\qbezier(22,46)(22,41)(25,39.5)

\qbezier(30,22)(26.5,21.7)(25,20.5)

\put(32,22){\line(-1,0){2}}
\qbezier(30,6)(26.5,6.5)(26.5,7)

\qbezier(25,20.5)(22,19)(22,14)
\qbezier(22,14)(22.7,9.3)(23.6,9)

\put(3,40){\vector(1,0){10}}
\put(47,40){\vector(-1,0){10}}
\put(0,43){\small $x$}
\put(46,13){\small $y$}

\put(3,20){\vector(1,0){10}}
\put(47,20){\vector(-1,0){10}}
\put(0,13){\small $x$}
\put(46,43){\small $y$}

\put(14,44){\small $*$}
\put(34,12){\small $*$}

\put(21,57){\scriptsize (1)}
\put(21,0){\scriptsize (2)}

\end{picture}
\end{minipage}
\hspace{-7pt}
\begin{minipage}{60pt}
\begin{picture}(50,60)
\put(6,6){\line(0,1){48}}
\put(6,6){\line(1,0){14}}
\put(6,54){\line(1,0){14}}
\put(18,22){\line(0,1){16}}

\qbezier(20,54)(23,53.7)(23.3,53.3)

\qbezier(26.5,51.5)(26.8,51)(27.5,49)
\qbezier(25,39.5)(27,39)(28,44)
\qbezier(28,44)(31,50)(33,46)
\qbezier(30,44.5)(33,42)(33,46)
\qbezier(27.5,49)(28,47.5)(28,46.5)

\put(18,38){\line(1,0){2}}
\qbezier(20,38)(23.5,38.5)(25,39.5)

\put(18,22){\line(1,0){2}}
\qbezier(20,22)(23.5,21.7)(25,20.5)

\qbezier(20,6)(23,6)(23,6.8)

\qbezier(25,20.5)(28,19)(28,14)
\qbezier(28,14)(28,9.5)(26,8.5)

\put(44,6){\line(0,1){48}}
\put(44,6){\line(-1,0){14}}
\put(44,54){\line(-1,0){14}}
\put(32,22){\line(0,1){16}}

\qbezier(30,22)(27,21.7)(26.7,21.3)

\qbezier(23.5,19.5)(23.2,19)(22.5,17)
\qbezier(25,7.5)(23,7)(22,12)
\qbezier(22,12)(19,18)(17,14)
\qbezier(20,12.5)(17,10)(17,14)
\qbezier(22.5,17)(22,15.5)(22,14.5)

\put(32,22){\line(-1,0){2}}
\qbezier(30,6)(26.5,6.5)(25,7.5)

\put(32,38){\line(-1,0){2}}
\qbezier(30,54)(26.5,53.7)(25,52.5)

\qbezier(30,38)(27,38)(27,38.8)

\qbezier(25,52.5)(22,51)(22,46)
\qbezier(22,46)(22,41.5)(24,40.5)

\put(3,30){\vector(1,0){10}}
\put(47,30){\vector(-1,0){10}}
\put(0,33){\small $x$}
\put(46,33){\small $y$}

\end{picture}
\end{minipage}

\hspace{-7pt}
\begin{minipage}{60pt}
\begin{picture}(50,60)
\put(6,6){\line(0,1){48}}
\put(6,6){\line(1,0){14}}
\put(6,54){\line(1,0){14}}
\put(18,22){\line(0,1){16}}

\qbezier(20,54)(23.5,53.7)(23.5,53.6)

\put(18,38){\line(1,0){2}}
\qbezier(20,38)(23.5,38.5)(25,39.5)

\qbezier(26.5,51.8)(28,51)(28,46)
\qbezier(28,46)(28,41)(25,39.5)

\put(18,22){\line(1,0){2}}
\qbezier(20,22)(23.5,21.7)(25,20.5)

\qbezier(20,6)(23.5,6.5)(23,6.7)

\qbezier(25,20,5)(28,19)(28,14)
\qbezier(28,14)(28,9)(26,8)

\put(44,6){\line(0,1){48}}
\put(44,6){\line(-1,0){14}}
\put(44,54){\line(-1,0){14}}
\put(32,22){\line(0,1){16}}

\qbezier(23.7,19.8)(22,19)(22,14)
\qbezier(22,14)(22,9)(25,7.5)

\qbezier(30,22)(27,21.7)(26.7,21.3)

\put(32,22){\line(-1,0){2}}
\qbezier(30,6)(26.5,6.5)(25,7.5)

\put(32,38){\line(-1,0){2}}
\qbezier(30,54)(26.5,53.7)(25,52.5)

\qbezier(30,38)(27,38)(27,38.8)

\qbezier(25,52.5)(22,51)(22,46)
\qbezier(22,46)(22,41.5)(24,40.5)

\end{picture}
\end{minipage}
\hspace{-10pt}
\begin{minipage}{40pt}
\begin{picture}(40,60)
\put(6,6){\line(0,1){48}}

\qbezier(6,6)(18,10)(18,30)
\qbezier(6,54)(18,50)(18,30)

\put(34,6){\line(0,1){48}}

\qbezier(34,6)(22,10)(22,30)
\qbezier(34,54)(22,50)(22,30)

\put(3,30){\vector(1,0){10}}
\put(37,30){\vector(-1,0){10}}
\put(0,33){\small $x$}
\put(36,33){\small $y$}

\end{picture}
\end{minipage}

\end{center}
\caption{a motion picture of $M$}
\label{Fig.J}
\end{figure}
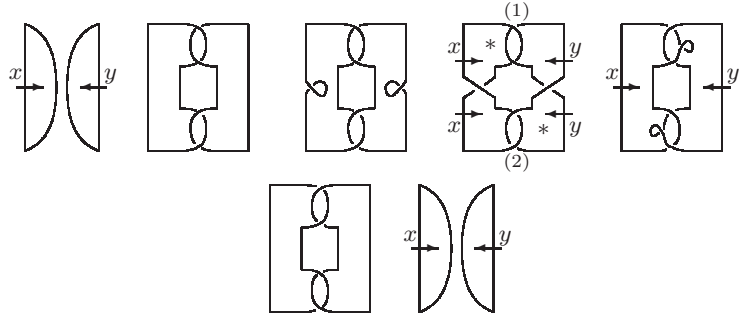

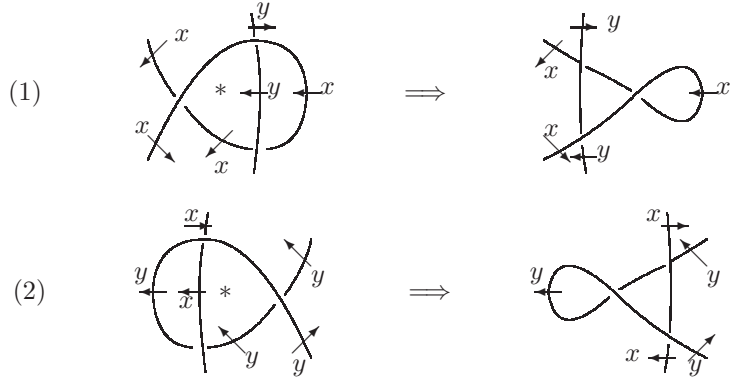
\begin{figure}[H]
\begin{center}
(1) \hspace{1cm}
\begin{minipage}{70pt}
\begin{picture}(70,60)

\qbezier(5,5)(25,50)(45,50)
\qbezier(45,50)(64,50)(65,30)
\qbezier(65,30)(64,10)(50,9)
\qbezier(44,9)(30,10)(19,25)
\qbezier(15,29)(7,40)(5,50)
\qbezier(44,60)(45,55)(45,52)
\qbezier(46,48)(49,30)(45,0)

\put(30,29){$*$}

\put(5,14){\vector(1,-1){10}}
\put(12,50){\vector(-1,-1){10}}
\put(37,16){\vector(-1,-1){10}}
\put(70,30){\vector(-1,0){10}}
\put(50,30){\vector(-1,0){10}}
\put(43,55){\vector(1,0){10}}

\put(0,15){$x$}
\put(15,50){$x$}
\put(30,0){$x$}
\put(70,30){$x$}
\put(50,30){$y$}
\put(46,60){$y$}

\end{picture}
\end{minipage}
\hspace{1cm}$\Longrightarrow $\hspace{1cm}
\begin{minipage}{70pt}
\begin{picture}(70,60)

\qbezier(5,5)(25,15)(40,30)
\qbezier(40,30)(62,50)(65,30)
\qbezier(65,30)(62,10)(42,28)
\qbezier(38,31)(30,36)(21,40)
\qbezier(18,42)(8,48)(5,50)
\qbezier(20,60)(18,40)(19,15)
\qbezier(19.5,10)(20,5)(21,0)

\put(5,14){\vector(1,-1){10}}
\put(12,50){\vector(-1,-1){10}}
\put(70,30){\vector(-1,0){10}}
\put(25,6){\vector(-1,0){10}}
\put(15,55){\vector(1,0){10}}

\put(5,14){$x$}
\put(5,36){$x$}
\put(70,30){$x$}
\put(25,6){$y$}
\put(29,55){$y$}

\end{picture}
\end{minipage}
\vspace{0.5cm}

(2) \hspace{1cm}
\begin{minipage}{70pt}
\begin{picture}(70,60)

\qbezier(65,5)(45,50)(25,50)
\qbezier(25,50)(6,50)(5,30)
\qbezier(5,30)(6,10)(20,9)
\qbezier(26,9)(40,10)(51,25)
\qbezier(55,29)(63,40)(65,50)
\qbezier(26,60)(25,55)(25,52)
\qbezier(24,48)(21,30)(25,0)

\put(30,28){$*$}

\put(58,7){\vector(1,1){10}}
\put(65,40){\vector(-1,1){10}}
\put(40,7){\vector(-1,1){10}}
\put(10,30){\vector(-1,0){10}}
\put(25,30){\vector(-1,0){10}}
\put(17,55){\vector(1,0){10}}

\put(58,0){$y$}
\put(65,34){$y$}
\put(40,3){$y$}
\put(-2,35){$y$}
\put(15,24){$x$}
\put(17,57){$x$}

\end{picture}
\end{minipage}
\hspace{1cm}$\Longrightarrow $\hspace{1cm}
\begin{minipage}{70pt}
\begin{picture}(70,60)

\qbezier(65,5)(45,15)(30,30)
\qbezier(30,30)(8,50)(5,30)
\qbezier(5,30)(8,10)(28,28)
\qbezier(32,31)(40,36)(49,40)
\qbezier(52,42)(62,48)(65,50)
\qbezier(50,60)(52,40)(51,15)
\qbezier(50.5,10)(50,5)(49,0)

\put(58,4){\vector(1,1){10}}
\put(65,40){\vector(-1,1){10}}
\put(10,30){\vector(-1,0){10}}
\put(53,5){\vector(-1,0){10}}
\put(48,55){\vector(1,0){10}}

\put(58,0){$y$}
\put(65,34){$y$}
\put(-2,35){$y$}
\put(34,4){$x$}
\put(42,57){$x$}

\end{picture}
\end{minipage}

\end{center}

\caption{Reidemeister moves of type III}
\label{Fig.K}
\end{figure}

Let $M= M_1 \cup M_2$ be a $2$-component compact surface properly embedded in ${\mathbb R}^3 \times [0,1]$ whose motion picture is given in Figure~\ref{Fig.J}.  The first still shows a  trivial $2$-component link, which is equal to the last still.  
The second still is obtained from the previous still by a pair of Reidemeister moves of type II.  The third is obtained by a pair of Reidemeister moves of type I.  There are two saddle points between the third still and the forth.  The fifth still is obtained from the forth by a pair of Reidemeister moves of type III, which are applied around the asterisked regions as in Figure~\ref{Fig.K}. 
The sixth still is obtained by a pair of Reidemeister moves of type I, and the last still is obtained by a pair of Reidemeister moves of type II.


For a positive integer $n$, put  
$$M^{(n)} = \{ x + (0,0,0,k) \in {\mathbb R}^4 ={\mathbb R}^3 \times  {\mathbb R} \, | \, 
x \in M, k \in \{0,1, \dots, n-1\} \}, $$
which is a properly embedded surface in ${\mathbb R}^3 \times [0, n]$.  Let $F=F^{(n)}$ be a closed surface in ${\mathbb R}^4 = {\mathbb R}^3 \times {\mathbb R}$ which is obtained from $M^{(n)}$ by attaching some disks trivially in $ {\mathbb R}^3 \times (-\infty, 0]$ and $
 {\mathbb R}^3 \times [n, \infty)$.   Then $F$ is a $2$-component surface-link, and each component is a non-orientable surface.  Satoh \cite{Satoh2002} proved that the minimal triple point number of $F$ is $2n$ by a geometric argument on the normal Euler number.   We can generalize this fact by the quandle cocycle invariant given in Example~\ref{example:cocycledihedral4}.  

\begin{theorem}\label{theorem:triple}
Let $F=F^{(n)}$ be the surface-link described above.  For any orientable surface-knot $K$, $t(F \# K) \geq 2n$.  Moreover, if $K$ is an orientable pseudo-ribbon surface-knot, then $t(F \# K) = 2n$.  {\rm  (A surface-knot $K$ is said to be {\it pseudo-ribbon} if $t(K)=0$.)   } 
\end{theorem}

\begin{proof} 
Considering  the motion picture of $M$ (Figure~\ref{Fig.J}) to be 
a $1$-parameter family of classical link diagram, we have a  broken surface diagram of $M$ with two triple points which correspond to the Reidemeister moves of type III in Figure~\ref{Fig.K}.  
Combining $n$ copies of the broken surface diagram and capping off by trivial disks, we obtain a diagram of the surface-link $F$ with $2n$ triple points, say $D$. 
  Let $(X, \rho)$, $Y$, and $\theta$ be as  in Example~\ref{example:cocycledihedral4}.    Let $C$ be a coloring of $D$ whose restriction to each copy of  the diagram of $M$ is as  in Figure~\ref{Fig.J},  where $x$ and $y$ are elements of $X$.  Then $c_{D,C} = n(e, x,y,x) - n(e,y,x,y) \in C_3^{ {\rm Q}, \rho}(X)_Y$.   If $(x,y) = (e_1, e_2)$, then  
$\theta (c_{D,C}) = 2n$.  By Proposition~\ref{proposition:cocycledihedral4}, $t(F) \geq 2n$.  
 Thus we have $t(F)= 2n$. 

It is seen by Proposition~\ref{proposition:cocycledihedral4} that $t (F \# K) \geq 2n$ where $F \# K$ is a connected sum of  $F$ and an orientable surface-knot $K$, since the colored diagram $(D,C)$ of $F$ can be extended uniquely to a colored diagram $(D', C')$ of $F \# K$ such that $c_{D,C} = c_{D', C'} \in C_3^{ {\rm Q}, \rho}(X)_Y$.   

If $K$ is an orientable pseudo-ribbon surface-knot, then $F\# K$ has a diagram with $2n$ triple points.  Thus $t(F \# K) = 2n$.  \end{proof}

The following  was suggested by S. Satoh. 

\begin{theorem}\label{theorem:P2P2}
For any positive integers $n_1$ and $n_2$ with $n_1 \equiv n_2 \, ({\rm mod}~2)$,  there is a surface-link $F_1 \cup F_2$ such that $F_i$ $(i=1,2)$ is the connected sum of $n_i$ copies of projective planes and $t(F_1 \cup F_2)= 2 \min\{n_1, n_2\}$.  
\end{theorem}

\begin{proof} 
Let $n= \min\{n_1, n_2\}$  
and $F= F^{(n)}$ the surface-link constructed above.  Let $K$ be an unknotted orientable surface-knot of genus $(\max\{n_1, n_2\} - n) / 2$.   As shown above, $t(F \#K)= 2n$.  The components of $F \# K$ are non-orientable surfaces which are the connected sums of $n_i$ $(i=1,2)$ copies of projective planes.  
\end{proof}\\

 
\noindent 
{\bf Acknowledgments}

The authors would like to thank Shin Satoh for helpful suggestions. 
This paper was completed while the second author visited Scott Carter and Masahico Saito.
The authors would also like to thank them for their helpful comments.


\begin{thebibliography}{99}
\bibitem{AG03} N. Andruskiewitsch and M. Gra\~na, {\it From racks to pointed Hopf algebras}, Adv. Math. {\bf 178} (2003), 177--243. 

\bibitem{AS05} S. Asami and S. Satoh, {\it An infinite family of non-invertible surfaces in $4$-space}, Bull. London Math. Soc. {\bf 37} (2005), 285--296. 

\bibitem{CEGS05} J. S. Carter, M. Elhamdadi, M. Gra\~na and M. Saito, {\it Cocycle knot invariants from quandle modules and generalized quandle cohomology}, Osaka J. Math. {\bf 42} (2005), 499--541. 

\bibitem{CJKLS03} J. S. Carter, D. Jelsovsky, S. Kamada, L. Langford and M. Saito, {\it Quandle cohomology and state-sum invariants of knotted curves and surfaces}, Trans. Amer. Math. Soc. {\bf 355} (2003), 3947--3989.

\bibitem{CJKS01b} J. S. Carter, D. Jelsovsky, S. Kamada and M. Saito, {\it Quandle homology groups, their Betti numbers, and virtual knots}, J. Pure Appl. Algebra {\bf 157} (2001), 135--155.

\bibitem{CKS01} J. S. Carter,  S. Kamada and M. Saito, 
{\it Geometric interpretations of quandle homology and cocycle knot invariants}, J. Knot Theory Ramifications 
{\bf 10} (2001) 345--358. 


\bibitem{CKSS01} J. S. Carter,  S. Kamada, M. Saito and S. Satoh, 
{\it A theorem of Sanderson on link bordisms in dimension $4$}, 
Algebr. Geom. Topol. {\bf 1} (2001) 299--310. 


\bibitem{CarterSaito} J. S. Carter and M. Saito, {\it Knotted surfaces and their diagrams}, Mathematical Surveys and Monographs, 55, Amer. Math. Soc., Providence, RI, 1998. 

\bibitem{FR92} R. Fenn and C. Rourke, {\it Racks and links in codimension two}, J. Knot Theory Ramifications {\bf 1} (1992), 343--406.  

\bibitem{FRS95} R. Fenn, C. Rourke and B. Sanderson, {\it Trunks and classifying spaces}, Appl. Categ. Structures {\bf 3} (1995), 321--356. 


\bibitem{FRS07} R. Fenn, C. Rourke and B. Sanderson, {\it The rack space}, Trans. Amer. Math. Soc. {\bf 359} (2007), 701--740.

\bibitem{I} M. Iwakiri, {\it Calculation of dihedral quandle cocycle invariants of twist spun 2-bridge knots}, J. Knot Theory Ramifications {\bf 14} (2005), 217--229.

\bibitem{Joyce} D. Joyce, {\it A classifying invariants of knots, the knot quandle}, J. Pure Appl. Algebra {\bf 23} (1982), 37--65.

\bibitem{Kamada2001} 
S. Kamada, 
{\it Wirtinger presentations for higher dimensional manifold knots obtained from diagrams},  
Fund. Math.  {\bf 168} (2001), 105--112. 


\bibitem{ka} S. Kamada, {\it Quandles with good involutions, their homologies and knot invariants}, in: Intelligence of Low Dimensional Topology 2006, Eds. J. S. Carter et. al., pp. 101-- 108, World Scientific Publishing Co., 2007. 

\bibitem{LN03} R. A. Litherland and S. Nelson, {\it The Betti numbers of some finite racks}, J. Pure Appl. Algebra {\bf 178} (2003), 187--202.


\bibitem{Matveev} S. Matveev, {\it Distributive groupoids in knot theory} (Russian), 
Mat. Sb. (N.S.) {\bf 119} (1982), 78--88; English translation: Math. USSR-Sb. {\bf 47} (1984), 73--83.

\bibitem{M} T. Mochizuki, {\it Some calculations of cohomology groups of finite Alexander quandles}, J. Pure Appl. Algebra {\bf 179} (2003), 287--330.

\bibitem{Roseman} D. Roseman, {\it Reidemeister-type moves for surfaces in four dimensional space}, in: Knot Theory, Banach Center Publications {\bf 42} (1998), 347--380. 

\bibitem{RS1} C. Rourke and B. Sanderson, {\it There are two 2-twist spun trefoils}, preprint, 
arxiv:math.GT/0006062:v1.

\bibitem{RS} C. Rourke and B. Sanderson, {\it A new classification of links and some calculation using it}, preprint, arxiv:math.GT/0006062:v2.

\bibitem{Oshiro} K. Oshiro, {\it Homology groups of trivial quandles with good involutions and triple linking numbers of surface-links}, preprint. 

\bibitem{Satoh2002} S. Satoh, {\it Triple point invariants of non-orientable surface-links}, Topology Appl. {\bf 121} (2002) 207--218. 

\bibitem{SatohA} S. Satoh, {\it A note on the shadow cocycle invariant of a knot with a base point}, J. Knot Theory Ramifications {\bf 16} (2007), 959--967.  


\bibitem{SatohShima2004} S. Satoh and A. Shima, {\it The 2-twist-spun trefoil has the triple point number four}, Trans. Amer. Math. {\bf 356} (2004), 1007--1024.

\bibitem{SatohShima2005} S. Satoh and A. Shima, {\it Triple point numbers and quandle cocycle invariants of knotted surfaces in 4-space}, New Zealand J. Math.  {\bf 34}  (2005), 71--79. 

\bibitem{Takasaki} M. Takasaki, {\it Abstraction of symmetric transformations} (Japanese), Tohoku Math. J. {\bf 49} (1943), 145--207.

\end{thebibliography}
\end{document}